\documentclass[11pt]{amsart}

\usepackage[margin=1.15in]{geometry}
\usepackage{amsmath,amssymb,amsthm}
\usepackage{aliascnt}
\usepackage{mathrsfs}
\usepackage[colorlinks=true,linkcolor=blue,citecolor=blue,urlcolor=blue]{hyperref}
\usepackage[nameinlink,capitalize,noabbrev]{cleveref}

\makeatletter
\renewcommand{\subsection}{%
  \@startsection{subsection}{2}{\z@}%
    {.5\linespacing\@plus.7\linespacing}%
    {.3\linespacing}%
    {\normalfont\bfseries}}

\renewcommand{\paragraph}{%
  \@startsection{paragraph}{4}{\z@}%
    {.5\linespacing\@plus.2\linespacing}%
    {-.5em}%
    {\normalfont\bfseries}}

\let\toda@setaddresses\@setaddresses
\renewcommand{\@setaddresses}{%
  \par\bigskip\noindent
  \begin{minipage}{\textwidth}
    \toda@setaddresses
  \end{minipage}\par}
\makeatother

\newcommand{\C}{\mathbb C}
\newcommand{\R}{\mathbb R}
\newcommand{\Pro}{\mathbb P}
\newcommand{\D}{\mathbb D}
\newcommand{\dd}{\mathrm d}
\newcommand{\ii}{\sqrt{-1}}
\newcommand{\FS}{\mathrm{FS}}
\newcommand{\ee}{\mathsf e}
\newcommand{\ff}{\mathsf f}
\newcommand{\hh}{\mathsf h}

\newtheorem{theorem}{Theorem}[section]
\newaliascnt{proposition}{theorem}
\newtheorem{proposition}[proposition]{Proposition}
\aliascntresetthe{proposition}
\newaliascnt{lemma}{theorem}
\newtheorem{lemma}[lemma]{Lemma}
\aliascntresetthe{lemma}
\newaliascnt{corollary}{theorem}
\newtheorem{corollary}[corollary]{Corollary}
\aliascntresetthe{corollary}
\theoremstyle{definition}
\newaliascnt{definition}{theorem}
\newtheorem{definition}[definition]{Definition}
\aliascntresetthe{definition}
\theoremstyle{remark}
\newaliascnt{remark}{theorem}
\newtheorem{remark}[remark]{Remark}
\aliascntresetthe{remark}
\numberwithin{equation}{section}

\title{Conical Singularities of Toda Solutions with One Finite-Energy Component and Classification on the Punctured Plane}
\author{Fangyuan Li}\author{Jingyu Mu}\author{Zhaohu Nie}\author{Yiqian Shi}\author{Chunhui Wei}\author{Bin Xu}

\makeatletter
\g@addto@macro\addresses{\setlength{\bigskipamount}{3pt}}
\makeatother
\address[Fangyuan Li, Yiqian Shi, and Bin Xu]{School of Mathematical Sciences, University of Science and Technology of China, Hefei 230026, China}
\address[Jingyu Mu]{Jintong Huayuan, No. 72 Xuesong Street, Aihui District, Heihe, Heilongjiang 164300, China}
\address[Zhaohu Nie]{Department of Mathematics and Statistics, Utah State University, Logan, UT 84322-3900, USA}
\address[Chunhui Wei]{School of Mathematics, Zhejiang University, Hangzhou 310058, China}

\makeatletter
\g@addto@macro\addresses{\par\addvspace{8pt}}
\makeatother
\email{%
  \normalfont
  \href{mailto:squarecircle31415@mail.ustc.edu.cn}{\nolinkurl{squarecircle31415@mail.ustc.edu.cn}} (F.~Li),
  \href{mailto:zhaohu.nie@usu.edu}{\nolinkurl{zhaohu.nie@usu.edu}} (Z.~Nie),
  \href{mailto:yqshi@ustc.edu.cn}{\nolinkurl{yqshi@ustc.edu.cn}} (Y.~Shi),
  \href{mailto:chunhuiwei@zju.edu.cn}{\nolinkurl{chunhuiwei@zju.edu.cn}} (C.~Wei),
  \href{mailto:bxu@ustc.edu.cn}{\nolinkurl{bxu@ustc.edu.cn}} (B.~Xu).
}


\begin{document}

\begin{abstract}
We study Toda systems associated with arbitrary complex simple Lie algebras. For solutions on a punctured disc, we prove that finite energy of at least one component near the puncture forces every component to have a conical singularity there.
On the punctured plane, we give a geometric formulation and proof of the Lin--Nie--Wei classification of solutions with prescribed cone exponents at the origin and finite total energy in at least one component. We express the parameter space as the complex centralizer of the monodromy modulo a maximal compact subgroup.
We also recover the formulas of Karmakar--Lin--Nie--Wei for the cone exponents at infinity and the total energies of the components.
\end{abstract}

\maketitle
\pagestyle{plain}

\section{Introduction}

Let $\mathfrak g$ be a complex simple Lie algebra of rank $n$, and let $\mathsf C=(c_{ij})$ be its Cartan matrix. We study the $\mathfrak g$-Toda system on a Riemann surface $S$, written globally as
\[
    \operatorname{Ric}(\boldsymbol\varphi)=2\mathsf C\boldsymbol\varphi.
\]
Here $\boldsymbol\varphi=(\varphi_1,\ldots,\varphi_n)^{\mathsf T}$ is a column vector of K\"ahler forms on $S$, and $\operatorname{Ric}(\boldsymbol\varphi)$ denotes the column vector of Ricci forms of these $\varphi_i$.

In a local holomorphic coordinate $z=x+\ii y$ on $S$, this system takes the form
\begin{equation}\label{eq:toda-local}
    \partial_z\partial_{\bar z}u_i
      =-\sum_{j=1}^n c_{ij}e^{u_j},
    \qquad
    \varphi_i=\frac{\ii}{2}e^{u_i}\dd z\wedge\dd\bar z,
    \quad 1\leq i\leq n.
\end{equation}
Here the $u_i$ are smooth real-valued functions. The energy of the $i$th component on a domain $U$ is its area $\int_U\varphi_i$, equal to $\int_U e^{u_i}\dd x\dd y$ in this coordinate.

\subsection{Conical singularities from finite energy of one component}

\begin{definition}\label{def:conical-singularity}
Let $\varphi$ be a K\"ahler form on a punctured neighborhood of a point $p$ on a Riemann surface. We say that $p$ is a \emph{conical singularity} of $\varphi$ with \emph{cone exponent} $\delta>0$ if there is a holomorphic coordinate $z$ centered at $p$ such that
\[
    \varphi=C|z|^{2\delta-2}(1+o(1))\frac{\ii}{2}\dd z\wedge\dd\bar z
    \qquad (z\to0)
\]
for some $C>0$. The corresponding metric has cone angle $2\pi\delta$. We allow $\delta=1$, so regular points are included. In particular, $\varphi$ has finite energy near $p$.
\end{definition}

For Toda solutions, we ask the converse question: does finite energy near an isolated puncture imply conical singularities?

In rank one, the Toda system is the Liouville equation for a metric of constant positive curvature. On a punctured disc, if the corresponding K\"ahler form has finite energy near the origin, Bryant \cite[Proposition~4]{Bryant} derives an explicit local expression for it, showing that the origin is a conical singularity.

Mu--Shi--Sun--Xu \cite[Theorem~1.2]{MuShiSunXu} describe solutions of the type $A_n$ Toda system with prescribed singular sources and finite energy in every component in terms of canonical unitary curves and derive explicit local expressions for the corresponding K\"ahler forms.
Mu \cite[Theorem~1.1 and Section~3.1]{MuThesis} considers solutions of the type $A_n$ Toda system on a punctured disc without prescribed singular sources. He proves that finite energy of one component near the puncture implies that every component has a conical singularity there.

The type $B_n$, $C_n$, and $G_2$ Toda systems arise as reductions of the type $A_{2n}$, $A_{2n-1}$, and $A_6$ Toda systems, respectively; see \cite[Lemmas~4.1 and~4.2]{NieBC} and \cite[(2.9)--(2.10)]{AoLinWei}. Consequently, Mu's conclusion also holds for these three types. These reductions do not cover the type $D_n$, $E_6$, $E_7$, $E_8$, and $F_4$ Toda systems; see \cite[Remark~1.11]{NieBC}.

Our first main result generalizes Mu's conclusion on conical singularities \cite[Theorem~1.1]{MuThesis} to Toda systems associated with every complex simple Lie algebra.

\begin{theorem}\label{thm:one-component-conical}
Let $\boldsymbol\varphi$ be a solution of the $\mathfrak g$-Toda system on the punctured unit disc $\D^*=\D\setminus\{0\}$. If one component of $\boldsymbol\varphi$ has finite energy near the origin, then every component has a conical singularity at the origin.
\end{theorem}

\subsection{Classification on the punctured plane}\label{subsec:intro-classification}

Our second question concerns the classification of solutions on the punctured plane $\C^*=\C\setminus\{0\}$. With prescribed singular sources at the origin and finite total energy in every component, the classification was obtained for the type $A_n$ Toda system by Lin--Wei--Ye \cite{LinWeiYe}, for types $B_n$ and $C_n$ by Nie \cite[Theorem~1.7]{NieBC}, and for type $G_2$ by Ao--Lin--Wei \cite[Theorem~2.1]{AoLinWei}. Mu \cite[Theorem~1.2 and Section~3.4]{MuThesis} classifies solutions of the type $A_n$ Toda system assuming that at least one component has finite total energy.

For Toda systems on $\C^*$ associated with arbitrary complex simple Lie algebras, the preprint by Lin--Nie--Wei \cite[Theorem~1.6]{LinNieWeiPreprint} classifies solutions with prescribed singular sources at the origin and finite total energy in every component. Its published version, coauthored with Karmakar, gives a representation-theoretic formula for these solutions \cite[Theorem~1.5]{KarmakarLinNieWei} and uses it to determine their total energies and asymptotic behavior at infinity \cite[Theorems~1.1 and~1.2]{KarmakarLinNieWei}. It leaves the complete classification for future work.

We give a geometric proof of the Lin--Nie--Wei classification and recover the Karmakar--Lin--Nie--Wei formulas for the total energies and the cone exponents at infinity.

We fix the global coordinate $z=x+\ii y$ on $\C^*$ and write
\[
    \mathrm{En}_i(u):=\int_{\C^*}e^{u_i}\dd x\dd y
\]
for the total energy of the $i$th component of a solution $u=(u_i)$ of \eqref{eq:toda-local}. We write $\delta_{p,i}(u)$ for the cone exponent of the corresponding K\"ahler form $\varphi_i$ at $p\in\{0,\infty\}$, whenever it exists. At infinity, the exponent is computed in the coordinate $\zeta=1/z$.

\begin{definition}\label{def:solution-spaces}
For $\delta=(\delta_1,\ldots,\delta_n)\in(\R_{>0})^n$, let $\mathscr S_\delta$ be the set of solutions $u$ of \eqref{eq:toda-local} on $\C^*$ satisfying $\delta_{0,i}(u)=\delta_i$ for every $i$ and $\mathrm{En}_j(u)<\infty$ for at least one $j$. Let $\mathscr R_\delta\subset\mathscr S_\delta$ consist of the \emph{radial solutions}, namely those for which each $u_i(z)$ depends only on $|z|$.
\end{definition}

Applying \cref{thm:one-component-conical} at both ends, we see that every solution $u$ on $\C^*$ with at least one component of finite total energy belongs to $\mathscr S_\delta$ for some $\delta\in(\R_{>0})^n$ and that each of its components has a conical singularity at $\infty$.

We introduce here the Lie-theoretic notation needed below. Let $G$ be the connected complex adjoint group of $\mathfrak g$, with a maximal compact subgroup $K$ and a maximal torus $T\subset K$. Set $\mathfrak t=\operatorname{Lie}(T)$ and $\mathfrak h=\mathfrak t\otimes_\R\C$, and let $T_\C\subset G$ be the complexification of $T$.

Choose simple roots $\alpha_1,\ldots,\alpha_n\in\mathfrak h^*$ with coroots $\hh_1,\ldots,\hh_n\in\ii\mathfrak t$ and fundamental weights $\omega_1,\cdots,\omega_n\in\mathfrak h^*$. Set $c_{ij}=\alpha_i(\hh_j)$ for the Cartan matrix. Then $\alpha_i=\sum_jc_{ij}\omega_j$.

Write $\chi_i:T_\C\to\C^*$ for the root character corresponding to $\alpha_i$. Since $G$ is adjoint, the character lattice of $T_\C$ is the root lattice \cite[Corollary~21.8]{MilneAlgebraicGroups}. The simple root characters therefore give an algebraic group isomorphism
\begin{equation}\label{eq:root-character-isomorphism}
    T_\C\xrightarrow{\sim}(\C^*)^n,
    \qquad t\longmapsto(\chi_i(t))_{i=1}^n.
\end{equation}
Taking absolute values identifies $T\backslash T_\C$ with $(\R_{>0})^n$ via $[t]\mapsto(|\chi_i(t)|)_{i=1}^n$.

Write $\mathcal C_+:=\{A\in\ii\mathfrak t: \alpha_i(A)>0\text{ for all }i\}$ for the open positive Weyl chamber. Let $\kappa$ be the longest Weyl group element relative to the chosen simple roots. Then $-\kappa\mathcal C_+=\mathcal C_+$. We use $\langle\cdot,\cdot\rangle$ for the natural pairing between $\mathfrak h^*$ and $\mathfrak h$.

For $\delta\in(\R_{>0})^n$, let $A_\delta\in\mathcal C_+$ be determined by $\alpha_i(A_\delta)=\delta_i$ for $1\leq i\leq n$, and set $M_\delta:=\exp(2\pi\ii A_\delta)\in T$. For $H=G,K$, write $Z_H(M_\delta)$ for the centralizer of $M_\delta$ in $H$.

\begin{theorem}\label{thm:intro-classification}
Fix $\delta\in(\R_{>0})^n$ and the notation above.
\begin{enumerate}
    \item There is a bijection
    \begin{equation}\label{eq:intro-classification}
        \Phi_\delta:Z_K(M_\delta)\backslash Z_G(M_\delta)
          \longrightarrow\mathscr S_\delta.
    \end{equation}
    Its restriction to the subset $T\backslash T_\C$ induced by $T_\C\subset Z_G(M_\delta)$ is a bijection onto $\mathscr R_\delta$.
    \item For every $u\in\mathscr S_\delta$, the cone exponents at infinity are
    \[
        \delta_{\infty,i}(u)=\alpha_i(-\kappa A_\delta),
        \qquad 1\leq i\leq n.
    \]
    The total energies of the components are
    \[
        \mathrm{En}_i(u)=\pi\langle\omega_i-\kappa\omega_i,A_\delta\rangle,
        \qquad 1\leq i\leq n.
    \]
\end{enumerate}
\end{theorem}

This is a geometric reformulation of the Lin--Nie--Wei classification \cite[Theorem~1.6]{LinNieWeiPreprint}. Every solution in $\mathscr S_\delta$ is obtained by applying $c\in Z_G(M_\delta)$ to the developing map of a fixed radial solution with monodromy $M_\delta$; see \eqref{eq:finite-energy-classification}. Two parameters $c,c'\in Z_G(M_\delta)$ give the same solution exactly when $c'=kc$ for some $k\in Z_K(M_\delta)$, giving the bijection $\Phi_\delta$ from the right coset space. Equivalently, replacing $c$ with $c^{-1}$ gives a parametrization by $Z_G(M_\delta)/Z_K(M_\delta)$.

\subsection{Method and organization}

In \cref{sec:energy-propagation}, we propagate finite energy between adjacent components along the connected Dynkin diagram using circular means and Jensen's inequality.
The Toda--flag correspondence recalled in \cref{sec:toda-flag-correspondence} associates a nondegenerate unitary integral curve in the flag variety with each Toda solution. The K\"ahler forms $\varphi_i$ are recovered as pullbacks of Fubini--Study forms along the curve's projective components.
The local problem therefore reduces to studying unitary projective curves on $\D^*$ with finite area near the origin.
In \cref{sec:projective-conical-singularities}, we obtain single-valued reductions of these curves by removing suitable monodromy factors. We use the auxiliary Fubini--Study potentials to obtain a zero-counting estimate from finite area. We then extend the single-valued reductions holomorphically across the puncture.
Computing the leading terms of the pullback Fubini--Study forms shows that each K\"ahler form $\varphi_i$ has a conical singularity at the origin. This proves \cref{thm:one-component-conical}.

In \cref{sec:local-flag-reductions}, we remove suitable monodromy factors from the associated flag curves on $\C^*$ and use the preceding extension result to extend the resulting single-valued reductions holomorphically to the Riemann sphere.
In \cref{sec:radial-solutions}, we characterize radial solutions as those whose single-valued reductions can be chosen constant, and then parameterize them.
In \cref{sec:global-classification}, we construct centralizer deformations of radial solutions. We use normalized algebraic lifts and the Drinfeld--Sokolov normal form to obtain a curve model in terms of constant Lie-theoretic data. Comparing cone exponents and monodromy then shows that every solution in $\mathscr S_\delta$ is a centralizer deformation of a fixed radial solution. We finally compute the cone exponents at infinity and the total energies of the components.
This completes the proof of \cref{thm:intro-classification}.

\section{Propagation of finite local energy}\label{sec:energy-propagation}

For a solution of the $\mathfrak g$-Toda system on the punctured unit disc $\D^*=\D\setminus\{0\}$, finite energy of one component near the origin implies finite energy of every component there. The analytic core of this argument is independent of the Lie-theoretic structure, so we first establish a propagation lemma for a broader class of coefficient matrices. Applying it to the Cartan matrix and iterating along the connected Dynkin diagram of $\mathfrak g$ then yields the result for the Toda system.

\begin{lemma}\label{lem:one-step-energy-propagation}
Let $B=(b_{ij})$ be a real $n\times n$ matrix, and fix an index $i\in\{1,\ldots,n\}$. Assume that $b_{ij}\leq0$ for every $j\ne i$, and let $u_1,\ldots,u_n$ be smooth real-valued functions on $\D^*$ satisfying
\begin{equation}\label{eq:general-exponential-equation}
    \partial_z\partial_{\bar z}u_i=-\sum_{j=1}^n b_{ij}e^{u_j}.
\end{equation}
If, for some $r_0\in(0,1)$,
\[
    \int_{0<|z|<r_0}e^{u_i}\frac{\ii}{2}\dd z\wedge\dd\bar z<\infty,
\]
then
\[
    \int_{0<|z|<r_0}e^{u_j}\frac{\ii}{2}\dd z\wedge\dd\bar z<\infty
\]
for every $j$ such that $b_{ij}<0$.
\end{lemma}

\begin{proof}
Use logarithmic polar coordinates $z=e^{-R+\ii\theta}$ and put $R_0=-\log r_0$. Thus $R\to\infty$ corresponds to $z\to0$, and the Euclidean area measure takes the form $e^{-2R}\dd R\,\dd\theta$.

Define the circular means
\[
    m_k(R)=\frac1{2\pi}\int_0^{2\pi}u_k(e^{-R+\ii\theta})\dd\theta, \qquad
    a_k(R)=\frac1{2\pi}\int_0^{2\pi}e^{u_k(e^{-R+\ii\theta})}\dd\theta.
\]
Define the normalized energies
\[
    J_k(R)=\int_{R_0}^R e^{-2s}a_k(s)\dd s
    =\frac1{2\pi}\int_{\{e^{-R}<|z|<r_0\}}
      e^{u_k}\frac{\ii}{2}\dd z\wedge\dd\bar z.
\]
Therefore, the $k$-th component has finite energy on $0<|z|<r_0$ if and only if
\[
    J_k(\infty):=\lim_{R\to\infty}J_k(R)<\infty.
\]

Taking the circular mean of \eqref{eq:general-exponential-equation} over $|z|=e^{-R}$ and using $\Delta=4\partial_z\partial_{\bar z}$, we obtain
\[
    m_i''(R)=-4b_{ii}e^{-2R}a_i(R)+4\sum_{j\ne i}(-b_{ij})e^{-2R}a_j(R).
\]
Integrating over the finite interval $[R_0,R]$ yields
\begin{equation}\label{eq:mean-identity}
    m_i'(R)=m_i'(R_0)-4b_{ii}J_i(R)+4\sum_{j\ne i}(-b_{ij})J_j(R).
\end{equation}
All terms in the last sum are nonnegative. Suppose that for some $j$, we have $b_{ij}<0$ and $J_j(\infty)=\infty$. By the hypothesis, $J_i(\infty)<\infty$, so the diagonal term in \eqref{eq:mean-identity} is bounded. Meanwhile the $j$-th term tends to $+\infty$, hence $m_i'(R)\to+\infty$. In particular, there exists $R_1>R_0$ such that $m_i(R)\ge3R$ for all $R>R_1$.

For each $R$, Jensen's inequality for the normalized circular measure gives $a_i(R)\geq e^{m_i(R)}$. It follows that
\[
    J_i(\infty)\ge\int_{R_1}^{\infty}e^{-2R}e^{m_i(R)}\dd R
    \ge \int_{R_1}^{\infty}e^R\dd R=\infty,
\]
a contradiction. Thus every component $u_j$ with $b_{ij}<0$ has finite energy.
\end{proof}

\begin{theorem}\label{thm:energy-propagation}
Let $u$ be a smooth real solution of \eqref{eq:toda-local} on $\D^*$. If, for some $i\in\{1,\ldots,n\}$ and some $r_0\in(0,1)$,
\[
    \int_{0<|z|<r_0}e^{u_i}\frac{\ii}{2}\dd z\wedge\dd\bar z<\infty,
\]
then for every $j\in\{1,\ldots,n\}$ we have
\[
    \int_{0<|z|<r_0}e^{u_j}\frac{\ii}{2}\dd z\wedge\dd\bar z<\infty.
\]
\end{theorem}
\begin{proof}
Apply \cref{lem:one-step-energy-propagation} with $B=\mathsf C$. The off-diagonal entries of a Cartan matrix are nonpositive. Moreover, $c_{ij}<0$ if and only if $c_{ji}<0$, and the Dynkin diagram of a complex simple Lie algebra is connected. Iterating the lemma along paths in the Dynkin diagram proves that finite energy of any one component implies finite energy of every component.
\end{proof}

\begin{remark}\label{rem:energy-propagation-scope}
More generally, suppose that $B$ has nonpositive off-diagonal entries and that $u_1,\ldots,u_n$ are smooth real-valued functions on $\D^*$ satisfying \eqref{eq:general-exponential-equation} for every $i$. Consider the directed graph with an arrow $i\to j$ whenever $b_{ij}<0$. Iterating \cref{lem:one-step-energy-propagation} shows that finite local energy of $u_i$ implies finite local energy of every component reachable from $i$. In particular, if this graph is strongly connected, finite local energy of any one component implies finite local energy of all components.

For semisimple $\mathfrak g$, the Cartan matrix is block diagonal, with one block for each simple factor. Applying \cref{thm:energy-propagation} to each block gives finite local energy of all components, provided at least one component in each block has finite local energy.

Even when the directed graph associated with $B$ is strongly connected, finite local energy need not imply conical singularities. For example, take
\[
    B=\begin{pmatrix}2&-3\\-3&2\end{pmatrix},
    \qquad
    u_1=u_2=\log\frac1{2|z|^2(\log|z|)^2}.
\]
They satisfy
\[
    \partial_z\partial_{\bar z}u_i=e^{u_i}
    =-\sum_j b_{ij}e^{u_j},
    \qquad i=1,2.
\]
Moreover, both components have finite energy near the origin:
\[
    \int_{0<|z|<r_0}
    e^{u_i}\frac{\ii}{2}\dd z\wedge\dd\bar z
    =\frac{\pi}{-\log r_0}<\infty.
\]
Nevertheless, the associated forms have cusp singularities at the origin, rather than conical singularities. The Toda--flag correspondence and \cref{thm:finite-area-conical} provide the additional input needed to obtain conical singularities for Toda solutions.
\end{remark}
\section{Toda systems and integral curves}\label{sec:toda-flag-correspondence}
We use the classical Toda--flag correspondence recalled by Shi--Wei--Xu \cite[Section~1]{WeiXuReduced} and record the parts needed below. We retain the Lie-theoretic notation introduced in \cref{subsec:intro-classification}. In \cite{WeiXuReduced}, our $G,K$ are denoted by $G_{\mathrm{ad}},K_{\mathrm{ad}}$, and our Cartan-matrix entries $c_{ij}$ by $a_{ji}$.

Let $\Delta^+$ be the positive root system determined by the chosen simple roots, and set
\[
    \mathfrak n_\pm:=\bigoplus_{\alpha\in\Delta^+}\mathfrak g_{\pm\alpha},
    \qquad N_\pm:=\exp(\mathfrak n_\pm).
\]
Let $B_0$ be the negative Borel subgroup given by
\[
    B_0=T_{\C}\ltimes N_-,\quad
    \mathfrak b_0=\mathfrak h\oplus\mathfrak n_-,\quad
    K\cap B_0=T.
\]

Let $F=G/B_0$ be the flag variety of Borel subgroups of $G$, with base point $B_0$. Let $\widetilde G\to G$ be the simply connected covering group, and let $\widetilde K$ be the inverse image of $K$. Our $\widetilde G,\widetilde K,B_0$ correspond to $G,K,B_{\mathrm{ad}}$ in \cite{WeiXuReduced}. With our negative Borel convention, let $V_i$ be the fundamental representation of $\widetilde G$ with lowest weight $-\omega_i$. We denote the action on $V_i$ and its differential by
\[
    \varrho_i:\widetilde G\longrightarrow\mathrm{GL}(V_i),\qquad
    \dd\varrho_i:\mathfrak g\longrightarrow\mathfrak{gl}(V_i).
\]

The finite center of $\widetilde G$ acts on each $V_i$ by scalars, so projectivization descends to $\overline\varrho_i:G\to\mathrm{PGL}(V_i)$. Let $\ell_i\in\Pro(V_i)$ be the lowest-weight line. The subgroup $B_0$ stabilizes this line \cite[Theorem~22.18]{MilneAlgebraicGroups}, giving the $i$th fundamental projection
\[
    \pi_i:F\longrightarrow\Pro(V_i),\qquad
    \pi_i(gB_0)=\overline\varrho_i(g)\cdot\ell_i.
\]
The projections $\pi_i$ are $G$-equivariant. Any element fixing all the $\ell_i$ belongs to $B_0$, since the Weyl group element indexing its double coset in the Bruhat decomposition \cite[Theorem~21.80]{MilneAlgebraicGroups} fixes every $-\omega_i$ and is therefore the identity. Thus the product of the projections
\begin{equation}\label{eq:fundamental-product-embedding}
    \Pi:F\longrightarrow\prod_{i=1}^n\Pro(V_i),\qquad
    x\longmapsto\bigl(\pi_1(x),\ldots,\pi_n(x)\bigr)
\end{equation}
identifies $F$ with the $G$-orbit of $(\ell_1,\ldots,\ell_n)$, which is closed because $F$ is projective.

Fix a $\widetilde K$-invariant Hermitian inner product on each $V_i$, unique up to a positive scalar, and denote the induced norms by $\|\cdot\|$. Then $\dd\varrho_i(A)$ is Hermitian for every $A\in\ii\mathfrak t$. Let $\varphi_i^{\FS}$ be the induced Fubini--Study form on $\Pro(V_i)$. Rescaling the inner product does not change $\varphi_i^{\FS}$. Denoting the natural projection by $[\ \cdot\ ]:V_i\setminus\{0\}\to\Pro(V_i)$, we normalize $\varphi_i^{\FS}$ by
\begin{equation}\label{eq:fs-normalization}
    [\ \cdot\ ]^*\varphi_i^{\FS}
    =\frac{\ii}{2}\partial\bar\partial\log\|\cdot\|^2.
\end{equation}
This is the normalization for which the Chern curvature of the induced metric on $\mathcal O_{\Pro(V_i)}(1)$ equals $-2\ii\varphi_i^{\FS}$. It is used for all Fubini--Study forms throughout the paper.

For $B\in F$, write $\mathfrak b=\operatorname{Lie}(B)$. Then $T_BF\simeq\mathfrak g/\mathfrak b$. The standard distribution $\mathcal N\subset TF$ is
\[
    \mathcal N(B)
    =\{X\in\mathfrak g:[X,[\mathfrak b,\mathfrak b]]
      \subset\mathfrak b\}/\mathfrak b.
\]
For each simple root $\alpha_j$, let $\mathfrak p_j(B)$ be the Lie algebra of the corresponding minimal parabolic subgroup containing $B$, and set $\mathcal N_j(B):=\mathfrak p_j(B)/\mathfrak b$. Then
\[
    \mathcal N=\bigoplus_{j=1}^n\mathcal N_j, \qquad
    \mathcal N_j(B_0)=\mathfrak g_{\alpha_j}\bmod\mathfrak b_0.
\]
We write $\operatorname{pr}_j:\mathcal N\to\mathcal N_j$ for the associated projection.

We now introduce the nondegenerate integral curves used in the Toda--flag correspondence. Throughout this section, let $S$ be a connected Riemann surface, and $p:\widetilde S\to S$ be the universal covering map of $S$.

\begin{definition}\label{def:nondegenerate-integral-curve}
A holomorphic flag curve $f:S\to F$ is \emph{integral} if $\dd f(TS)\subset f^*\mathcal N$. An integral curve $f$ is \emph{nondegenerate} if, for every $j\in\{1,\ldots,n\}$, the simple-root component $\operatorname{pr}_j\circ\ \dd f$ is nowhere vanishing.
\end{definition}

\begin{definition}\label{def:unitary-flag-curve}
A \emph{unitary flag curve} on $S$ is a pair $f=(\widetilde f,\rho_f)$ consisting of a holomorphic map $\widetilde f:\widetilde S\to F$ and a homomorphism $\rho_f:\pi_1(S)\to K$ such that
\[
    \widetilde f(\gamma\cdot\widetilde z)
    =\rho_f(\gamma)\cdot\widetilde f(\widetilde z)
\]
for every $\gamma\in\pi_1(S)$ and $\widetilde z\in\widetilde S$. We call $\widetilde f$ the \emph{developing map} of $f$ and $\rho_f$ its \emph{monodromy representation}.

A unitary flag curve is called \emph{integral} if its developing map is integral. A unitary integral curve is called \emph{nondegenerate} if its developing map is nondegenerate.
\end{definition}

The following lemma implies that the developing map of a nondegenerate unitary integral curve determines its monodromy representation.

\begin{lemma}[{\cite[Lemma~2.8]{WeiXuReduced}}]\label{lem:pointwise-stabilizer}
Let $f:S\to F$ be a nondegenerate integral curve. For $c\in G$, if $c$ fixes $f$ pointwise, then $c$ is the identity.
\end{lemma}

We use the following form of the classical Toda--flag correspondence \cite[Corollary~1.10]{WeiXuReduced}. The formula \eqref{eq:toda-component-forms} follows from Positsel'skii's generalized Pl\"ucker formula \cite{Positselskii}.

\begin{theorem}[Toda--flag correspondence]
\label{thm:toda-flag-correspondence}
Let $\operatorname{Int}_G(S)_K$ be the space of nondegenerate unitary integral curves on $S$, and let $\operatorname{Toda}_{\mathfrak g}(S)$ be the space of $\mathfrak g$-Toda solutions on $S$. There is a bijection
\[
    K\backslash\operatorname{Int}_G(S)_K
    \xrightarrow{\ \sim\ }
    \operatorname{Toda}_{\mathfrak g}(S).
\]
The left $K$-action is given by
\[
    k\cdot(\widetilde f,\rho_f)
    =\bigl(k\cdot\widetilde f,k\rho_f k^{-1}\bigr).
\]
The class of $f=(\widetilde f,\rho_f)$ is sent to the Toda solution $\boldsymbol\varphi=(\varphi_1,\ldots,\varphi_n)$ determined by
\begin{equation}\label{eq:toda-component-forms}
    p^*\varphi_i
    =(\pi_i\circ\widetilde f)^*\varphi_i^{\FS},
    \qquad 1\leq i\leq n.
\end{equation}
\end{theorem}

Since the singularity analysis treats the fundamental projective components one at a time, we introduce the following notion for a general finite-dimensional Hermitian vector space $V$.

\begin{definition}\label{def:unitary-projective-curve}
A \emph{unitary projective curve} on $S$ is a pair $f=(\widetilde f,\rho_f)$ consisting of a holomorphic map $\widetilde f:\widetilde S\to\Pro(V)$ and a homomorphism $\rho_f:\pi_1(S)\to\mathrm{PU}(V)$ such that
\[
    \widetilde f(\gamma\cdot\widetilde z)
       =\rho_f(\gamma)\cdot\widetilde f(\widetilde z)
\]
for every $\gamma\in\pi_1(S)$ and $\widetilde z\in\widetilde S$. We call $\widetilde f$ its \emph{developing map} and $\rho_f$ its \emph{monodromy representation}.
\end{definition}

For a unitary projective curve $f=(\widetilde f,\rho_f)$,  the $\mathrm{PU}(V)$-invariance of $\varphi^{\mathrm{FS}}$ implies that the form $\widetilde f^{\,*}\varphi^{\FS}$ descends to a well-defined nonnegative $(1,1)$-form on $S$, denoted by $f^*\varphi^{\FS}$. It may vanish at branch points.

For a unitary flag curve $f=(\widetilde f,\rho_f)$, its $i$-th \emph{fundamental projective component} is
\[
    f_i=(\widetilde f_i,\rho_{f_i}),\qquad
    \widetilde f_i:=\pi_i\circ\widetilde f,\qquad
    \rho_{f_i}:=\overline\varrho_i\circ\rho_f.
\]
The $G$-equivariance of $\pi_i$ and the $\widetilde K$-invariance of the Hermitian inner product on $V_i$ ensure that $f_i$ is a unitary projective curve. If $f$ corresponds to a Toda solution $\boldsymbol\varphi$,
then \eqref{eq:toda-component-forms} identifies $f_i^*\varphi_i^{\FS}$ with the K\"ahler form $\varphi_i$. Consequently, each developing map $\widetilde f_i$ is an immersion, and the local energy of each Toda component is the area of its fundamental projective component.

\section{Holomorphic extension and conical singularities}\label{sec:projective-conical-singularities}

We prove \cref{thm:one-component-conical} by characterizing nonconstant unitary projective curves on the punctured disc $\D^*=\D\setminus\{0\}$ with finite area near the origin.

\begin{theorem}\label{thm:finite-area-conical}
Let $V$ be a finite-dimensional Hermitian vector space, and let $f=(\widetilde f,\rho_f)$ be a unitary projective curve on $\D^*$ with nonconstant developing map $\widetilde f$. Then $f^*\varphi^{\FS}$ has finite area near the origin if and only if the origin is a conical singularity of $f^*\varphi^{\FS}$.
\end{theorem}

We fix the universal covering map of $\D^*$:
\begin{equation}\label{eq:universal-covering}
    p:\mathbb H\longrightarrow\D^*\ ,\quad w\longmapsto e^{w},
\end{equation}
where $\mathbb H$ is the \emph{left} half-plane. Fix the positively oriented generator $\gamma_0\in\pi_1(\D^*)$, acting on $\mathbb H$ by $\gamma_0\cdot w=w+2\pi\ii$. For a unitary projective curve $f=(\widetilde f,\rho_f)$ on $\D^*$, we call
\[
    M_f:=\rho_f(\gamma_0)\in \mathrm{PU}(V)
\]
the \emph{monodromy} of the curve.

To analyze the behavior of $f$ near the puncture, we reduce to a holomorphic extension problem for a single-valued curve on $\D^*$ by removing the monodromy. Under the finite-area hypothesis, we prove that this curve extends across $0$ and use a holomorphic lift of the extension to compute the leading terms of $f^*\varphi^{\FS}$.

\begin{definition}\label{def:expansive-monodromy-factor}
Let $M\in\mathrm U(V)$. An \emph{expansive monodromy factor associated with $M$} is a holomorphic map $\phi:\mathbb H\to\mathrm{GL}(V) $ such that
\begin{equation}\label{eq:expansive-monodromy-factor-norm}
    \phi(\gamma_0\cdot w)=M\phi(w),\qquad
    \|\phi(w)v\|\geq\|v\|,
\end{equation}
for every $w\in\mathbb H$ and $v\in V$.
\end{definition}

\begin{lemma}\label{lem:expansive-monodromy-factor}
Every $M\in\mathrm U(V)$ admits an expansive monodromy factor. More precisely, there is a Hermitian endomorphism $B\in\operatorname{End}(V)$ such that
\[
    M=\exp(2\pi\ii B)\ , \qquad
    \operatorname{Spec}(B)\subset(-1,0],
\]
and then
\[
    \phi_B(w):=\exp\bigl(Bw\bigr)
\]
is an expansive monodromy factor associated with $M$. We call $\phi_B$ the \emph{logarithmic expansive monodromy factor associated with $B$}.
\end{lemma}
\begin{proof}
By the spectral theorem, choose an orthonormal eigenbasis $(e_j)$ of $M$ and write $Me_j=e^{2\pi\ii b_j}e_j$, where $b_j\in(-1,0]$. Define $B$ by $Be_j=b_je_j$. Then $B$ is Hermitian, has the stated spectrum, and $M=\exp(2\pi\ii B)$. Hence
\[
    \phi_B(w+2\pi\ii)=e^{2\pi\ii B}\phi_B(w)=M\phi_B(w).
\]
In the basis $(e_j)$, write $v=\sum_jv_je_j$. For $z=e^w$, we have
\[
    \|\phi_B(w)v\|^2=\sum_j |z|^{2b_j}|v_j|^2
    \geq\sum_j |v_j|^2=\|v\|^2,
\]
because $0<|z|<1$ and $b_j\leq0$. This proves the required norm inequality.
\end{proof}

\begin{definition}\label{def:contractive-single-valued-reduction}
Let $f=(\widetilde f,\rho_f)$ be a unitary projective curve. Choose a lift $\widehat M_f\in\mathrm U(V)$ of the monodromy $M_f$, and fix an expansive monodromy factor $\phi$ associated with $\widehat M_f$. The map
\[
    w\longmapsto \phi(w)^{-1}\cdot\widetilde f(w)
\]
is invariant under the deck transformations and therefore descends to a holomorphic curve $h:\D^*\to\Pro(V)$. We call $h$ the \emph{contractive single-valued reduction} of $f$ relative to $\phi$.
\end{definition}

The following proposition generalizes Mu's extension result \cite[Proposition~3.1]{MuThesis} by replacing his diagonal power factors with expansive monodromy factors on $V$.

\begin{proposition}\label{thm:holomorphic-extension}
Let $V$ be a finite-dimensional Hermitian vector space, and let $f=(\widetilde f,\rho_f)$ be a unitary projective curve on $\D^*$. Fix a lift $\widehat M_f\in\mathrm U(V)$ of its monodromy and an expansive monodromy factor $\phi$ associated with $\widehat M_f$, and let $h:\D^*\to\Pro(V)$ be the corresponding contractive single-valued reduction. If, for some $r\in(0,1)$,
\begin{equation}\label{eq:finite-area}
    \int_{0<|z|<r}f^*\varphi^{\FS}<\infty,
\end{equation}
then $h$ extends holomorphically to a map $\overline h:\D\to\Pro(V)$.
\end{proposition}

Following Mu's approach, we introduce the auxiliary Fubini--Study potentials $q_\ell$, inspired by Wu \cite[Theorem~4.1]{Wu}, which relate the Fubini--Study area of $f$ to the hyperplane intersections of $h$. We first construct compatible holomorphic lifts.

\begin{lemma}\label{lem:compatible-lifts}
Let $f=(\widetilde f,\rho_f)$, $\widehat M_f$, $\phi$, and $h$ be as in \cref{def:contractive-single-valued-reduction}. Then $\widetilde f$ and $h$ admit holomorphic lifts $\vec f:\mathbb H\to V\setminus\{0\}$ and $\vec h:\D^*\to V\setminus\{0\}$, respectively, such that
\[
    \vec f(\gamma_0\cdot w)=\widehat M_f\vec f(w)\ , \quad
    \vec f=\phi\cdot(\vec h\circ p).
\]
The norm of $\vec f$ descends to a positive smooth function $\|\vec f\|(z):=\|\vec f(w)\|$ for $z\in\D^*$ and any $w\in p^{-1}(z)$. On $\D^*$, we have the identity
\begin{equation}\label{eq:single-valued-norm}
    \ii\partial\bar\partial\log\|\vec f\|
    =f^*\varphi^{\FS}.
\end{equation}
Moreover,
\begin{equation}\label{eq:compatible-lifts-norm}
    \|\vec f\|(z)\geqslant\|\vec h(z)\|,
    \qquad z\in\D^*.
\end{equation}
\end{lemma}
\begin{proof}
The line bundle $h^*\mathcal O_{\Pro(V)}(-1)$ is holomorphically trivial because $\D^*$ is a non-compact Riemann surface \cite[Theorem~30.3]{Forster}. Choose a nowhere-vanishing holomorphic section of this bundle. Viewed as a section of $\D^*\times V$, it gives a holomorphic lift $\vec h$ of $h$. Set
\[
    \vec f(w):=\phi(w)\vec h(p(w)),
    \qquad[\vec f(w)]=\widetilde f(w).
\]
The equivariance of $\phi$ and the deck invariance of $\vec h\circ p$ give the equivariance of $\vec f$.

The unitary equivariance of $\vec f$ makes its norm deck invariant. For a local inverse $s:U\to\mathbb H$ of $p$, the descended norm equals $\|\vec f\circ s\|$, so it is positive and smooth. Applying \eqref{eq:fs-normalization} to this local lift gives
\[
    \ii\partial\bar\partial\log\|\vec f\|\big|_U
    =(\widetilde f\circ s)^*\varphi^{\FS}
    =f^*\varphi^{\FS}\big|_U,
\]
which proves \eqref{eq:single-valued-norm}. The norm inequality \eqref{eq:compatible-lifts-norm} follows from the norm condition in \eqref{eq:expansive-monodromy-factor-norm}, which explains the term \emph{contractive} in \cref{def:contractive-single-valued-reduction}.
\end{proof}

\begin{lemma}\label{lem:auxiliary-fs-potentials}
Let $f$ and the compatible lifts $\vec f,\vec h$ be as in \cref{lem:compatible-lifts}. Write $\|\cdot\|_*$ for the dual norm on $V^*$. For a fixed $\ell\in V^*\setminus\{0\}$ with $\ell(\vec h)\not\equiv0$, define
\[
    q_\ell(z)=\log\frac{\|\ell\|_*\|\vec f\|(z)}{|\ell(\vec h(z))|}.
\]
This function has the following properties:
    \begin{enumerate}
        \item $q_\ell(z)\geqslant 0$;
        \item Its singular set is discrete. Near each singular point $z_0\in\D^*$,
        \begin{equation}\label{eq:auxiliary-fs-singularity}
            q_\ell(z)=-m\log|z-z_0|+O(1),
        \end{equation}
        where $m$ is the zero multiplicity of $\ell(\vec h)$ at $z_0$ and the remainder is smooth. Away from the singular set,
        \begin{equation}\label{eq:auxiliary-fs-curvature}
            \frac12\dd\dd^cq_\ell=f^*\varphi^{\FS}\ ,
            \qquad \dd^c=\ii(\bar\partial-\partial);
        \end{equation}
        \item The circular mean
        \[
            Q_\ell(R)=\frac1{2\pi}\int_0^{2\pi}q_\ell(e^{-R+\ii\theta})\dd\theta\ ,
            \qquad R>0
        \]
        is locally absolutely continuous on $(0,\infty)$. Moreover, it is $C^1$ away from the locally finite set $\mathcal R_\ell=\bigl\{-\log|z|:\ell(\vec h(z))=0\bigr\}$.
    \end{enumerate}
\end{lemma}
\begin{proof}
By \eqref{eq:compatible-lifts-norm}, for $z\in\D^*$ we have
\[
    |\ell(\vec h(z))|
    \leq \|\ell\|_*\|\vec h(z)\|
    \leq \|\ell\|_*\|\vec f\|(z).
\]
This proves the nonnegativity of $q_\ell$. Since $\ell(\vec h(z))$ is holomorphic and not identically zero, its zero set is discrete. If $z_0$ is a zero of multiplicity $m$, then $\ell(\vec h(z))=(z-z_0)^m h_0(z)$ for some holomorphic function $h_0(z)$ with $h_0(z_0)\ne0$. Thus, near $z_0$,
\[
    q_\ell(z)=-m\log|z-z_0|+\psi(z), \qquad
    \psi(z):=\log\|\ell\|_*+\log\|\vec f\|(z)-\log|h_0(z)|.
\]
The remainder $\psi(z)$ extends smoothly across $z_0$, since $\|\vec f\|$ is positive and smooth by \cref{lem:compatible-lifts} and $h_0$ is nonvanishing near $z_0$. Away from these zeros, $\log|\ell(\vec h(z))|$ is harmonic. Since $\dd\dd^c=2\ii\partial\bar\partial$, the identity \eqref{eq:single-valued-norm} gives
\[
    \frac12\dd\dd^c q_\ell=\ii\partial\bar\partial\log\|\vec f\|
    =f^*\varphi^{\FS}.
\]

On a neighborhood of any compact annulus $K_{a,b}=\{e^{-b}\leq|z|\leq e^{-a}\}$, subtract the finitely many logarithmic singularities to write $q_\ell=q-\sum_\nu m_\nu\log|z-z_\nu|$ with $q$ smooth. Jensen's formula \cite[Chapter~XI, Theorem~1.2]{Conway} then gives, for $R\in[a,b]$,
\[
    Q_\ell(R)=Q(R)+\sum_\nu m_\nu\min\{R,R_\nu\},
    \qquad R_\nu=-\log|z_\nu|,
\]
where $Q$ is the smooth circular mean of $q$. The finite sum is absolutely continuous and is $C^1$ away from the $R_\nu$. Since $[a,b]$ is arbitrary, this proves the last assertion.
\end{proof}

\begin{proof}[Proof of \cref{thm:holomorphic-extension}]
\leavevmode\par
Choose the lifts $\vec f,\vec h$ in \cref{lem:compatible-lifts}. Choose $\ell_0\in V^*$ such that $\ell_0(\vec h)\not\equiv0$, and extend $\ell_0$ to a basis $\ell_0,\ldots,\ell_N$ of $V^*$, where $N=\dim V-1$. It is enough to prove that every quotient
\[
    v_j=\frac{\ell_j(\vec h)}{\ell_0(\vec h)}\ ,
    \qquad 1\leq j\leq N
\]
is meromorphic at $0$. Clearing poles and cancelling any common vanishing factor then gives a holomorphic extension of $h$.

Let $\ell\in V^*$ satisfy $\ell(\vec h)\not\equiv0$. We show that $\ell(\vec h)$ has only finitely many zeros near $0$. We use the auxiliary Fubini--Study potential $q_\ell$ and its circular mean $Q_\ell$ from \cref{lem:auxiliary-fs-potentials}. Choose $0<r_0<r$ so that $\ell(\vec h)$ has no zero on $|z|=r_0$, and put $R_0=-\log r_0$. For $R>R_0$, let $n_\ell(R)$ be the total multiplicity of the zeros of $\ell(\vec h)$ in the annulus $A_R=\{e^{-R}\leq|z|\leq r_0\}$.

For $R>R_0$ outside $\mathcal R_\ell$, delete disjoint discs $D_\varepsilon(z_\nu)$ around the zeros $z_\nu\in A_R$. With every circle below parametrized counterclockwise, Stokes' theorem gives
\[
\int_{A_R\setminus\bigcup_\nu D_\varepsilon(z_\nu)}\dd\dd^c q_\ell
=\int_{|z|=r_0}\dd^c q_\ell-\int_{|z|=e^{-R}}\dd^c q_\ell
-\sum_\nu\int_{\partial D_\varepsilon(z_\nu)}\dd^c q_\ell.
\]
The outer, inner, and deleted-circle boundary integrals are, respectively
\[
    \int_{|z|=r_0}\dd^c q_\ell=:2\pi C_\ell, \qquad
    -2\pi Q_\ell'(R), \qquad
    -2\pi m_\nu+o(1) \quad(\varepsilon\to0),
\]
where $C_\ell$ is independent of $R$ and $m_\nu$ is the
zero multiplicity at $z_\nu$. The last $o(1)$ estimate follows from \eqref{eq:auxiliary-fs-singularity}. Letting $\varepsilon\to0$ and using \eqref{eq:auxiliary-fs-curvature} yields
\begin{equation}\label{eq:zero-counting-identity}
    \frac1\pi\int_{e^{-R}\leq|z|<r_0}f^*\varphi^{\FS}-C_\ell
    =Q_\ell'(R)+n_\ell(R).
\end{equation}
Since $\mathcal R_\ell$ is locally finite, this identity holds for almost every $R>R_0$. By \eqref{eq:finite-area}, the area term in \eqref{eq:zero-counting-identity} is bounded above. Thus there is a constant $M_\ell\geq0$ such that
\[ Q_\ell'(R)+n_\ell(R)\leq M_\ell \]
for almost every $R>R_0$. For any $T>R_1>R_0$, the absolute continuity and nonnegativity of $Q_\ell$ and the monotonicity of $n_\ell$ give
\[
    n_\ell(R_1)(T-R_1)
    \leq \int_{R_1}^T n_\ell(R)\dd R
    \leq M_\ell(T-R_1)+Q_\ell(R_1).
\]
Letting $T\to\infty$ gives $n_\ell(R_1)\leq M_\ell$. Thus $\ell(\vec h)$ has only finitely many zeros near $0$.

Apply this conclusion first to $\ell=\ell_0$. After shrinking the punctured disc, $\ell_0(\vec h)$ has no zeros there, so every $v_j$ is holomorphic on the smaller punctured disc. If some nonconstant $v_j$ had an essential singularity at $0$, the great Picard theorem \cite[Chapter~XII, Theorem~4.2]{Conway} would give a finite value $\lambda$ assumed infinitely often near $0$. Since $v_j$ is nonconstant, $(\ell_j-\lambda\ell_0)(\vec h)$ is not identically zero, but it has infinitely many zeros near $0$, a contradiction. Thus every $v_j$ is meromorphic at $0$.
\end{proof}

\begin{proof}[Proof of \cref{thm:finite-area-conical}]
\leavevmode\par
Conical singularities directly imply finite area. Conversely, suppose that $f^*\varphi^{\FS}$ has finite area near the origin. Choose a unitary lift $\widehat M_f$ of the monodromy and its Hermitian logarithm $B$ as in \cref{lem:expansive-monodromy-factor}. Fix an orthonormal eigenbasis of $B$ and write
\[
    B=\operatorname{diag}(a_0,\ldots,a_N)\ ,\qquad
   a_j\in(-1,0]\ ,\ N=\dim V-1.
\]
Let $\phi_B(w)=e^{Bw}$ be the logarithmic expansive monodromy
factor, and let $h$ be the corresponding contractive single-valued reduction. By \cref{thm:holomorphic-extension}, $h$ extends holomorphically across the puncture. Let
\[
    \vec h=(h^0,\ldots,h^N)
\]
be a nowhere-vanishing holomorphic lift of this extension near $0$.

Choose a branch of $\log z$ on a simply connected open set $U\subset\D^*$ with $0\in\overline U$. By \cref{def:contractive-single-valued-reduction}, the local projective map $z\mapsto\widetilde f(\log z)$ has the holomorphic lift
\begin{equation}\label{eq:monodromy-adapted-local-lift}
    \vec f=(f^0,\ldots,f^N)
    =(z^{a_0}h^0,\ldots,z^{a_N}h^N).
\end{equation}
Changing the logarithm branch multiplies this lift by a constant unitary transformation, which does not affect the pullback $f^*\varphi^{\FS}$.

For every $h^j\not\equiv0$, write
\[
    h^j(z)=z^{m_j}(c_j+O(z))\ ,\qquad
    m_j\in\mathbb Z_{\geq0}\ ,\ c_j\ne0,
\]
and set
\[
    \lambda=\min_{h^j\not\equiv0}(a_j+m_j).
\]
The lowest-order terms in the squared norm give, for some $C_0>0$,
\begin{equation}\label{eq:norm-leading-order}
    \|\vec f(z)\|^2=C_0|z|^{2\lambda}(1+o(1))\qquad\text{as }z\to0.
\end{equation}

On $U$, the normalization of the Fubini--Study form \eqref{eq:fs-normalization} gives
\begin{equation}\label{eq:fs-wronskian-formula}
    f^*\varphi^{\FS}\big|_U
    =\frac{1}{\|\vec f\|^4}\sum_{0\leq j<k\leq N}
      |f^j(f^k)'-f^k(f^j)'|^2
      \frac{\ii}{2}\dd z\wedge\dd\bar z.
\end{equation}
Differentiating \eqref{eq:monodromy-adapted-local-lift}, for $j<k$ we obtain
\[
    f^j(f^k)'-f^k(f^j)'=z^{a_j+a_k-1}H_{jk}(z)\ ,\quad
    H_{jk}=(a_k-a_j)h^jh^k+z(h^j(h^k)'-h^k(h^j)').
\]
Each $H_{jk}$ is single-valued and holomorphic near $0$, and at least one is not identically zero. Otherwise $f^j(f^k)'-f^k(f^j)'\equiv0$ for all $j<k$. Every ratio $f^j/f^k$ would then be constant wherever $f^k\ne0$, and the identity theorem would make $\widetilde f$ constant, contrary to the hypothesis.

For every $H_{jk}\not\equiv0$, write
\[
    H_{jk}(z)=z^{\nu_{jk}}(d_{jk}+O(z))\ ,
    \qquad \nu_{jk}\in\mathbb Z_{\geq0}\ ,\ d_{jk}\ne0,
\]
and set
\[
    \delta=\min_{H_{jk}\not\equiv0}\bigl(a_j+a_k+\nu_{jk}-2\lambda\bigr).
\]
The lowest-order terms in the numerator of \eqref{eq:fs-wronskian-formula} give, for some $C_1>0$,
\begin{equation}\label{eq:wronskian-leading-order}
    \sum_{j<k}|f^j(f^k)'-f^k(f^j)'|^2
    =C_1|z|^{4\lambda+2\delta-2}(1+o(1)).
\end{equation}
Substituting \eqref{eq:norm-leading-order} and \eqref{eq:wronskian-leading-order} into \eqref{eq:fs-wronskian-formula} gives, with $C=C_1/C_0^2>0$,
\begin{equation}\label{eq:conical-asymptotic}
    f^*\varphi^{\FS}
    =C|z|^{2\delta-2}(1+o(1))\frac{\ii}{2}\dd z\wedge\dd\bar z
    \qquad (z\to0).
\end{equation}
Since $h^j$ and $H_{jk}$ are holomorphic near $0$ and $|z^{a_j}|$ is independent of $\arg z$, the asymptotic holds uniformly in $\arg z$ on $\D^*$. The finite-area hypothesis then forces $\delta>0$, so the origin is a conical singularity of $f^*\varphi^{\FS}$.
\end{proof}

\begin{proof}[Proof of \cref{thm:one-component-conical}]
By \cref{thm:energy-propagation}, every component has finite energy near the origin. Choose an associated nondegenerate unitary integral curve $f$ by \cref{thm:toda-flag-correspondence}. Its fundamental projective components $f_i$ satisfy $f_i^*\varphi_i^{\FS}=\varphi_i$ by \eqref{eq:toda-component-forms}. Thus each $f_i$ has finite area near the origin and a nonconstant developing map, since $\varphi_i$ is K\"ahler. Applying \cref{thm:finite-area-conical} to each $f_i$ proves the assertion.
\end{proof}

\begin{remark}
When $\dim V=2$, \cref{thm:holomorphic-extension} recovers the meromorphic extension step in Bryant's proof of \cite[Proposition~4]{Bryant}. There the developing map is written as $z^\alpha\psi(z)$ with $\psi$ single-valued, and finite area implies that $\psi$ is meromorphic at $0$.

In our notation and normalization, Bryant's result gives the following exact local normal form in a suitable coordinate $\zeta$ centered at the puncture:
\[
    f^*\varphi^{\FS}
    =\frac{\delta^2|\zeta|^{2\delta-2}}
           {(1+|\zeta|^{2\delta})^2}
      \frac{\ii}{2}\dd\zeta\wedge\dd\bar\zeta.
\]
In the original coordinate $z$, this normal form gives the asymptotic behavior \eqref{eq:conical-asymptotic}.
\end{remark}
\section{Global factorization of unitary integral curves}
\label{sec:local-flag-reductions}

For a Toda solution on $\C^*$ with one component of finite total energy, \cref{thm:global-flag-reduction} gives a global factorization of an associated nondegenerate unitary integral curve. The key step is to extend its flag single-valued reduction relative to a logarithmic $G$-valued monodromy factor across $0$ and $\infty$. We then use the cone exponents at the origin to choose the logarithm $A$.

We use the universal cover $p:\C\to\C^*$, $p(w)=e^w$ for $\C^*$. Its restriction to $\mathbb H$ is the covering map \eqref{eq:universal-covering}. We use the same positively oriented generator $\gamma_0$ for $\pi_1(\D^*)$ and $\pi_1(\C^*)$. For a unitary flag curve $f=(\widetilde f,\rho_f)$ on $\D^*$ or $\C^*$, we call
\[
    M_f:=\rho_f(\gamma_0)\in K
\]
the \emph{monodromy} of the curve. We identify $\C^*$ with $\Pro^1\setminus\{0,\infty\}$ using the affine coordinate $z$.

\begin{theorem}\label{thm:global-flag-reduction}
Let $u$ be a smooth Toda solution on $\C^*$ with at least one component of finite total energy, let $f=(\widetilde f,\rho_f)$ be an associated nondegenerate unitary integral curve, and let $M\in K$ be its monodromy around $0$. Then, for every $k\in K$ and $A\in\ii\mathfrak t$ satisfying $k^{-1}Mk=e^{2\pi\ii A}$, there is a unique holomorphic map $\overline g:\Pro^1\to F$ such that
\begin{equation}\label{eq:global-factorization}
    \widetilde f(w)=k e^{Aw}\cdot\overline g(e^w),
    \qquad w\in\C.
\end{equation}
\end{theorem}

The following definitions clarify the meanings of the two factors in \eqref{eq:global-factorization}. We give the definitions on $\D^*$; for $\C^*$, replace $\D^*$ by $\C^*$ and $\mathbb H$ by $\C$.

\begin{definition}\label{def:G-valued-monodromy-factor}
Let $M\in K$. A \emph{$G$-valued monodromy factor associated with $M$} is a holomorphic map $\Phi:\mathbb H\to G$ satisfying
\[
    \Phi(\gamma_0\cdot w)=M\Phi(w).
\]
In particular, by conjugacy of maximal tori and surjectivity of $\exp:\mathfrak t\to T$, choose $k\in K$ and $A\in\ii\mathfrak t$ such that $k^{-1}Mk=e^{2\pi\ii A}$. Then
\begin{equation}\label{eq:logarithmic-G-valued-monodromy-factor}
    \Phi_{k,A}(w):=ke^{Aw}
\end{equation}
is a monodromy factor associated with $M$, called a \emph{logarithmic $G$-valued monodromy factor}.
\end{definition}

\begin{definition}\label{def:flag-single-valued-reduction}
Let $f=(\widetilde f,\rho_f)$ be a unitary flag curve on $\D^*$, and fix a $G$-valued monodromy factor $\Phi$ associated with $M_f$. The map
\[
    w\longmapsto\Phi(w)^{-1}\cdot\widetilde f(w)
\]
is invariant under the deck transformations and therefore descends to a holomorphic curve $g:\D^*\to F$. We call $g$ the \emph{flag single-valued reduction} of $f$ relative to $\Phi$.
\end{definition}

Removing a monodromy factor does not in general preserve integrality. For $a\in G$, write $L_a(x):=a\cdot x$. For $X\in\mathfrak g$, define its infinitesimal action by
\[
    X_F(x):=\left.\frac{\dd}{\dd t}\right|_{t=0}
        \exp(tX)\cdot x.
\]

The left logarithmic derivative $\Phi^{-1}\dd\Phi$ is unchanged by constant left multiplication, so the equivariance of $\Phi$ makes it deck invariant. It therefore descends to a $\mathfrak g$-valued holomorphic one-form on $\D^*$, denoted by the same symbol. For a holomorphic map $g:\D^*\to F$, we write $(\Phi^{-1}\dd\Phi)_F(g)$ for the one-form induced by the infinitesimal action along $g$, namely
\[
    \bigl((\Phi^{-1}\dd\Phi)_F(g)\bigr)_z(v):=
    \bigl((\Phi^{-1}\dd\Phi)_z(v)\bigr)_F(g(z)),
    \qquad z\in\D^*,\quad v\in T_z\D^*.
\]

In the setting of \cref{def:flag-single-valued-reduction}, differentiating $\widetilde f=\Phi\cdot(g\circ p)$ gives
\begin{equation}\label{eq:twisted-differential}
    (L_{\Phi^{-1}})_*\dd\widetilde f
    =p^*\bigl(\dd g+(\Phi^{-1}\dd\Phi)_F(g)\bigr).
\end{equation}
Since $\mathcal N$ is $G$-invariant and $p$ is a surjective local biholomorphism, this shows that $f$ is integral if and only if
\[
    \dd g+(\Phi^{-1}\dd\Phi)_F(g)
    \quad\text{takes values in }g^*\mathcal N.
\]
The summands $\mathcal N_j$ are also $G$-invariant, and the projections $\operatorname{pr}_j$ are $G$-equivariant. When the integrality condition holds, the same identity shows that $f$ is nondegenerate if and only if every simple-root component of this one-form is nowhere vanishing.

To extend the flag single-valued reduction $g$, we compare its fundamental projective components with the contractive single-valued reductions in \cref{thm:holomorphic-extension}, which extend holomorphically under the finite-area hypothesis, and assemble the componentwise extensions.

\begin{lemma}\label{lem:laurent-comparison}
Let $f=(\widetilde f,\rho_f)$ be a unitary flag curve on $\D^*$. There are contractive single-valued reductions $h_i$ of the fundamental projective components $f_i$ such that, for every flag single-valued reduction $g$ of $f$ relative to a logarithmic $G$-valued monodromy factor and every $i\in\{1,\ldots,n\}$,
\[
    h_i\text{ extends holomorphically across }0
    \quad\Longleftrightarrow\quad
    \pi_i\circ g\text{ extends holomorphically across }0.
\]
\end{lemma}

\begin{proof}
For each $i$, choose a unitary lift $\widehat M_i\in\mathrm U(V_i)$ of the monodromy of $f_i$. Let $B_i$ and $\phi_{B_i}$ be obtained from $\widehat M_i$ by \cref{lem:expansive-monodromy-factor}, and let $h_i$ be the corresponding contractive single-valued reduction of $f_i$.

Let $g$ be any flag single-valued reduction relative to $\Phi_{k,A}$, and put $g_i:=\pi_i\circ g$ and $A_i:=\dd\varrho_i(A)$. Choose a lift $\widetilde k\in\widetilde K$ of $k$, and put $K_i:=\varrho_i(\widetilde k)\in\mathrm U(V_i)$. Since $K_i e^{2\pi\ii A_i}K_i^{-1}$ and $\widehat M_i$ are both unitary lifts of the monodromy of $f_i$, there is $s_i\in\R$ such that
\[
    K_i e^{2\pi\ii A_i}K_i^{-1}=e^{2\pi\ii s_i}\widehat M_i.
\]
Put $B_i^0:=K_i^{-1}B_iK_i$. Both $A_i$ and $B_i^0$ are Hermitian, and
\[
    e^{2\pi\ii B_i^0}
    =K_i^{-1}\widehat M_iK_i
    =e^{2\pi\ii(A_i-s_iI)}.
\]
Since $\operatorname{Spec}(B_i^0)\subset(-1,0]$ and $b\mapsto e^{2\pi\ii b}$ is injective on this interval, $B_i^0$ is scalar on each eigenspace of $e^{2\pi\ii(A_i-s_iI)}$. It therefore commutes with $A_i$, and the corresponding eigenvalues of $A_i-s_iI$ and $B_i^0$ differ by integers. Thus $N_i:=A_i-s_iI-B_i^0$ is Hermitian with integral spectrum. Set $\Phi_i(w):=K_i e^{A_iw}$, which induces on $\Pro(V_i)$ the action of $\Phi_{k,A}$. For $z=e^w$,
\[
    \phi_{B_i}(w)^{-1}\Phi_i(w)
      =K_ie^{-B_i^0w}e^{A_iw}
      =K_ie^{(N_i+s_iI)w}
      =e^{s_iw}K_i z^{N_i}.
\]
Hence a single-valued matrix $L_i(z):=K_i z^{N_i}$ represents the projective action of $\phi_{B_i}^{-1}\Phi_i$ on $\Pro(V_i)$. Since $N_i$ has integral spectrum, $L_i$ is Laurent.

Equivariance of $\pi_i$ gives $\Phi_i\cdot(g_i\circ p)=\widetilde f_i=\phi_{B_i}\cdot(h_i\circ p)$, so $h_i=L_i\cdot g_i$. Since both $L_i$ and $L_i^{-1}$ are Laurent matrices, applying either to a holomorphic lift and then clearing poles and cancelling any common vanishing factor gives the claimed equivalence of extensions.
\end{proof}

\begin{proof}[Proof of \cref{thm:global-flag-reduction}]
\leavevmode\par
Fix $k$ and $A$ as in the statement, and let $g:\C^*\to F$ be the flag single-valued reduction of $f$ relative to $\Phi_{k,A}$. By \cref{thm:energy-propagation} at both ends and \eqref{eq:toda-component-forms}, every fundamental projective component of $f$ has finite area near $0$ and $\infty$. For each $i$, choose $h_i$ as in \cref{lem:laurent-comparison}. By \cref{thm:holomorphic-extension}, the contractive single-valued reduction $h_i$ extends holomorphically across $0$, and \cref{lem:laurent-comparison} then gives a holomorphic extension $\overline g_i:\D\to\Pro(V_i)$ of $\pi_i\circ g$. The product of these extensions,
\[
    (\overline g_1,\ldots,\overline g_n):
    \D\longrightarrow\prod_{i=1}^n\Pro(V_i)
\]
agrees with $\Pi\circ g$ on $\D^*$, where $\Pi$ is the standard closed embedding \eqref{eq:fundamental-product-embedding}. Therefore these projective extensions assemble uniquely into a holomorphic extension $\overline g:\D\to F$ of $g$.

At infinity, put $\zeta=1/z=e^v$ and $v=-w$. The developing map $\widetilde f_\infty(v):=\widetilde f(-v)$ has monodromy $M^{-1}$ around $\zeta=0$. The logarithmic $G$-valued monodromy factor $\Phi_{k,-A}$ is associated with $M^{-1}$, and the corresponding flag single-valued reduction is
\[
    \Phi_{k,-A}(v)^{-1}\cdot\widetilde f_\infty(v)
    =\Phi_{k,A}(-v)^{-1}\cdot\widetilde f(-v)
    =g(e^{-v})=g(1/\zeta).
\]
The same argument in the coordinate $\zeta$, with logarithm $-A$, extends the same $g$ holomorphically across $\infty$. The two extensions uniquely give $\overline g:\Pro^1\to F$.
\end{proof}

By \cref{thm:one-component-conical}, every solution in \cref{thm:global-flag-reduction} belongs to $\mathscr S_\delta$ for some $\delta\in(\R_{>0})^n$. We now use these cone exponents to choose the logarithm $A$ in \eqref{eq:global-factorization}. We retain the notation $A_\delta,M_\delta$ from \cref{subsec:intro-classification}. Put $\mathfrak k=\operatorname{Lie}(K)$. Complete the fixed coroots $\hh_i$ to a Chevalley basis with
\[
    \ee_i\in\mathfrak g_{\alpha_i},\qquad
    \ff_i\in\mathfrak g_{-\alpha_i},\qquad
    [\ee_i,\ff_j]=\delta_{ij}\hh_i.
\]
The basis is chosen so that the Cartan involution $\sigma$ with fixed-point algebra $\mathfrak k$ satisfies $\sigma(\ee_i)=-\ff_i$, $\sigma(\ff_i)=-\ee_i$, and $\sigma(\hh_i)=-\hh_i$. The fundamental coweights are given by
\[
    E_i:=\sum_{j=1}^n c^{ji}\hh_j,\qquad
    \alpha_j(E_i)=\delta_{ij},\qquad (c^{ij}):=\mathsf C^{-1}.
\]
Given a Toda solution $u$, its Toda connection is given by (see \cite[\S1.3]{WeiXuReduced})
\begin{equation}\label{eq:toda-connection}
    \Theta=\frac12\sum_{i=1}^n E_i
       \bigl(\partial_z u_i\,\dd z-
             \partial_{\bar z}u_i\,\dd\bar z\bigr)
       +\sum_{i=1}^n e^{u_i/2}
          \bigl(\ee_i\,\dd z-\ff_i\,\dd\bar z\bigr).
\end{equation}
It is $\mathfrak k$-valued since $\sigma(\Theta)=\Theta$. The Toda equation is equivalent to its flatness \cite[Lemma~1.4]{WeiXuReduced}.

\begin{proposition}\label{prop:cone-adapted-logarithm}
In the setting of \cref{thm:global-flag-reduction}, if $u\in\mathscr S_\delta$, then $M$ is $K$-conjugate to $M_\delta$. In particular, one may take $A=A_\delta$ in \eqref{eq:global-factorization}.
\end{proposition}

\begin{proof}
Since $u\in\mathscr S_\delta$, in the global coordinate $z$ we have
\[
    u_i(z)=2(\delta_i-1)\log|z|+\log C_i+o(1)
    \qquad (z\to0).
\]
Set $v_i:=u_i-2(\delta_i-1)\log|z|-\log C_i$ and $\delta_*:=\min_j\delta_j>0$. On the fixed annulus $1/2<|\zeta|<2$, the cone asymptotics and \eqref{eq:toda-local} give $v_i(r\zeta)=o(1)$ and $\Delta_\zeta v_i(r\zeta)=O(r^{2\delta_*})$ uniformly as $r\to0$. The interior gradient estimate for Poisson's equation \cite[Theorem~3.9]{GilbargTrudinger} therefore yields $\sup_{|z|=r}r|\nabla v_i(z)|=o(1)$. Consequently
\begin{equation}\label{eq:cone-radial-derivative}
    r\partial_r u_i=2(\delta_i-1)+o(1)
\end{equation}
uniformly on $|z|=r$. We also have $|z|e^{u_i/2}=O(r^{\delta_i})\to0$ for $\delta_i>0$, so the angular component of \eqref{eq:toda-connection} along the positively oriented circle $z=re^{\ii\theta}$ satisfies
\begin{align*}
    \Theta(\partial_\theta)
      &=\frac{\ii}{2}\sum_{i=1}^n E_i
          \bigl(z\partial_z u_i+\bar z\partial_{\bar z}u_i\bigr)
        +\ii\sum_{i=1}^n e^{u_i/2}
          \bigl(z\ee_i+\bar z\ff_i\bigr)\\
      &=\ii(A_\delta-\rho^\vee)+o(1) \qquad(r\to 0)
\end{align*}
uniformly for $\theta\in[0,2\pi]$. Here $\rho^\vee\in\ii\mathfrak t$ is characterized by $\alpha_i(\rho^\vee)=1$ for $1\leq i\leq n$. Since $\chi_i(e^{2\pi\ii\rho^\vee})=1$ for every $i$, the injectivity of \eqref{eq:root-character-isomorphism} gives $e^{2\pi\ii\rho^\vee}=1$.

Let $U_r$ be the parallel transport along the parametrized circle, satisfying $U_r'=U_r\Theta(\partial_\theta)$ and $U_r(0)=1$. By the uniform estimate on $\Theta(\partial_\theta)$, we obtain
\begin{equation}\label{eq:cone-angular-transport}
    U_r(\theta)\longrightarrow
    e^{\ii\theta(A_\delta-\rho^\vee)}
    \qquad(r\to0),
\end{equation}
uniformly for $\theta\in[0,2\pi]$. Taking $\theta=2\pi$ therefore gives $U_r(2\pi)\to M_\delta$.

By the identification of monodromy with holonomy in \cite[proof of Theorem~1.7]{WeiXuReduced}, each $U_r(2\pi)$ is $K$-conjugate to $M$. Since conjugacy classes in $K$ are closed, the limit $M_\delta$ is also $K$-conjugate to $M$. Choose $k\in K$ with $k^{-1}Mk=M_\delta$. Applying \cref{thm:global-flag-reduction} with this $k$ and $A=A_\delta$ proves the final assertion.
\end{proof}
\section{Radial solutions and constant reductions}
\label{sec:radial-solutions}

We characterize the solutions in $\mathscr R_\delta$ defined in \cref{def:solution-spaces} by the existence of a constant flag single-valued reduction relative to a logarithmic monodromy factor, and then construct and parameterize these solutions. We retain the notation from \cref{subsec:intro-classification}.

\begin{theorem}\label{thm:radial-parameterization}
Fix $\delta\in(\R_{>0})^n$. A solution $u\in\mathscr S_\delta$ is radial if and only if $u=u^{A_\delta,x}$ for some constant reduction $x$ for $A_\delta$ (see \cref{def:constant-reduction}). Moreover, there is a bijection
\begin{equation}\label{eq:radial-parameterization}
    T\backslash T_\C\longrightarrow\mathscr R_\delta.
\end{equation}
\end{theorem}

We begin with the rotational equivariance of developing maps of nondegenerate unitary integral curves associated with radial solutions.

\begin{lemma}\label{lem:radial-rotation-lift}
Let $u$ be a radial Toda solution on $\C^*$ with cone exponents $\delta=(\delta_i)$ at $0$, and let $f=(\widetilde f,\rho_f)$ be an associated nondegenerate unitary integral curve with monodromy $M\in K$. Then there is $h\in K$ such that for all $w\in\C$ and $\theta\in[0,2\pi]$,
\[
    \widetilde f(w+\ii\theta)
    =h e^{\ii\theta A_\delta}h^{-1}\cdot\widetilde f(w),
    \qquad M =h M_\delta h^{-1}.
\]
\end{lemma}

\begin{proof}
Choose $\Psi:\C\to K$ satisfying $\Psi^{-1}\dd\Psi=p^*\Theta$, where $\Theta$ is the Toda connection \eqref{eq:toda-connection}. By the construction in \cite[proof of Theorem~1.7]{WeiXuReduced}, $\Psi B_0$ is the developing map of an associated curve. By \cref{thm:toda-flag-correspondence}, after a constant left $K$-translation of $\Psi$ we may assume $\Psi B_0=\widetilde f$.

For $\theta\in[0,2\pi]$, set $R_\theta(z)=e^{\ii\theta}z$ and $t_\theta=e^{\ii\theta\rho^\vee}\in T$. Radiality and $\alpha_i(\rho^\vee)=1$ give
\[
    R_\theta^*\Theta
    =\frac12\sum_i E_i
    \bigl(\partial_z u_i\dd z
            -\partial_{\bar z}u_i\dd\bar z\bigr)
    +\sum_i e^{u_i/2}
    \bigl(e^{\ii\theta}\ee_i\dd z
            -e^{-\ii\theta}\ff_i\dd\bar z\bigr)
    =\operatorname{Ad}_{t_\theta}\Theta.
\]
Hence $\Psi(w+\ii\theta)t_\theta$ and $\Psi(w)$ have the same left logarithmic derivative, and therefore differ by a constant left translation $k_\theta\in K$:
\begin{equation}\label{eq:radial-frame-equivariance}
    \Psi(w+\ii\theta)t_\theta=k_\theta\Psi(w).
\end{equation}
Since $t_\theta\in T\subset B_0$, projecting \eqref{eq:radial-frame-equivariance} to $F=G/B_0$ gives
\begin{equation}\label{eq:radial-flag-equivariance}
    \widetilde f(w+\ii\theta)=k_\theta\cdot\widetilde f(w).
\end{equation}

Let $\delta_*=\min_i\delta_i>0$. Along the positive real axis, radiality and the cone asymptotics at $0$ give
\[
    \Theta(\partial_r)=\sum_i e^{u_i(r)/2}(\ee_i-\ff_i)
       =O(r^{\delta_*-1}), \qquad(r\to0).
\]
It is integrable at $0$, so the equation $\frac{\dd}{\dd r}\Psi(\log r)=\Psi(\log r)\Theta(\partial_r)$ implies that $\Psi(\log r)\to h\in K$ as $r\to0$. Together with \eqref{eq:cone-angular-transport}, this gives, for each fixed $\theta\in[0,2\pi]$,
\[
    \Psi(\log r+\ii\theta)\longrightarrow he^{\ii\theta(A_\delta-\rho^\vee)}, \qquad(r\to0).
\]
Taking the limit in \eqref{eq:radial-frame-equivariance} identifies
\[
    k_\theta=he^{\ii\theta(A_\delta-\rho^\vee)}t_\theta h^{-1}
             =he^{\ii\theta A_\delta}h^{-1}.
\]
Combining this with \eqref{eq:radial-flag-equivariance} proves the asserted rotation formula.

Setting $\theta=2\pi$ in \eqref{eq:radial-flag-equivariance} and using \cref{lem:pointwise-stabilizer}, we obtain $M=k_{2\pi}=h M_\delta h^{-1}$.
\end{proof}

We now consider the simplest possible curves: $w\mapsto e^{A_\delta w}\cdot x$ for $x\in F$. Their flag single-valued reductions relative to $w\mapsto e^{A_\delta w}$ are the constant maps $g\equiv x$. It turns out that every solution in $\mathscr R_\delta$ arises from a curve of this form.

\begin{definition}\label{def:constant-reduction}
For $A\in\mathcal C_+$ and $x\in F$, let $f_{A,x}$ be the unitary flag curve with developing map $\widetilde f_{A,x}(w):=e^{Aw}\cdot x$ and monodromy $e^{2\pi\ii A}$. We call $x$ a \emph{constant reduction for $A$} if $f_{A,x}$ is a nondegenerate unitary integral curve. We write $u^{A,x}$ for the associated Toda solution.
\end{definition}

\begin{proposition}\label{prop:constant-reductions}
Let $A\in\mathcal C_+$ and $x\in F$. Then $x$ is a constant reduction for $A$ if and only if
\[
    A_F(x)\in\mathcal N(x), \qquad
    \operatorname{pr}_i A_F(x)\ne0
    \quad(1\leq i\leq n).
\]
In this case, $u^{A,x}$ is radial and $\delta_{0,i}(u^{A,x})=\alpha_i(A)$ for $1\leq i\leq n$.
\end{proposition}

\begin{proof}
For $\Phi(w)=e^{Aw}$ and $g\equiv x$, we have $\dd g+(\Phi^{-1}\dd\Phi)_F(g)=A_F(x)\dd z/z$. By \eqref{eq:twisted-differential}, the stated condition is equivalent to $x$ being a constant reduction for $A$. Rotation through $\theta$ gives
\begin{equation}\label{eq:constant-reduction-rotation}
    \widetilde f_{A,x}(w+\ii\theta)=e^{\ii\theta A}\cdot\widetilde f_{A,x}(w).
\end{equation}
Since $e^{\ii\theta A}\in K$, the Fubini--Study area forms of the fundamental projective components are rotation invariant. Hence $u^{A,x}$ is radial.

The local calculation in the proof of \cref{thm:finite-area-conical} shows that $u^{A,x}$ has conical singularities at both ends and hence has finite total energy. Let $\delta=(\delta_i)$ be its cone exponents at $0$. Applying \cref{lem:radial-rotation-lift} gives $h\in K$ such that
\[
    \widetilde f_{A,x}(w+\ii\theta)
      =he^{\ii\theta A_\delta}h^{-1}\cdot\widetilde f_{A,x}(w).
\]
Comparing this with \eqref{eq:constant-reduction-rotation} and applying \cref{lem:pointwise-stabilizer}, we obtain
\[
    e^{\ii\theta A}=he^{\ii\theta A_\delta}h^{-1}.
\]
Differentiating at $\theta=0$ yields $A=\operatorname{Ad}_hA_\delta$. Since $A,A_\delta\in\mathcal C_+$, we have $A=A_\delta$, which proves $\delta_{0,i}(u^{A,x})=\alpha_i(A)$.
\end{proof}

We can now construct the constant reductions for a given $A\in\mathcal C_+$.

\begin{proposition}\label{prop:constant-reduction-construction}
For every $A\in\mathcal C_+$ and $b=(b_1,\ldots,b_n)\in(\C^*)^n$, there is a unique $n_A(b)\in N_+$ satisfying
\begin{equation}\label{eq:seed-root-coordinates}
    \operatorname{Ad}_{n_A(b)^{-1}}A
        =A+\sum_{i=1}^n b_i\ee_i.
\end{equation}
The points $x_A(b):=n_A(b)\cdot B_0\in F$ are exactly the constant reductions for $A$. The action of $T_\C$ on this locus is free and transitive, with
\[
    t\cdot x_A(b)
        =x_A\bigl(\chi_1(t)b_1,\ldots,\chi_n(t)b_n\bigr),
    \qquad t\in T_\C.
\]
\end{proposition}
\begin{proof}
By \cref{prop:constant-reductions}, a point $x=a\cdot B_0$ is a constant reduction precisely when
\[
    \operatorname{Ad}_{a^{-1}}A
      \in\mathfrak b_0+\sum_i\C\ee_i
\]
and every $\ee_i$-coefficient is nonzero. By the Bruhat decomposition \cite[Theorem~21.80]{MilneAlgebraicGroups}, write $x=v\widehat\iota B_0$, where $v\in N_+$ and $\widehat\iota\in N_K(T)$ represents an element $\iota$ of the Weyl group. Take $a=v\widehat\iota$. Since $\operatorname{Ad}_{v^{-1}}A\in A+\mathfrak n_+$, the $\mathfrak g_{\alpha_i}$-component of $\operatorname{Ad}_{a^{-1}}A$ vanishes whenever $\iota\alpha_i\in-\Delta^+$. Nondegeneracy therefore forces $\iota\alpha_i\in\Delta^+$ for every $i$. Thus $\iota$ preserves $\Delta^+$, and hence $\iota=1$. Thus $x=vB_0$ with $v$ unique, and \eqref{eq:seed-root-coordinates} holds with $n_A(b)=v$ for some $b\in(\C^*)^n$.

Conversely, fix $b\in(\C^*)^n$ and seek $n_A(b)=\exp Y$. The height of a positive root $\alpha=\sum_i d_i\alpha_i$
is $\sum_i d_i$. We grade $Y=\sum_{k\ge1}Y_k$ by this height. The height-$k$ component of $e^{-\operatorname{ad}Y}A$ is
\[
    [A,Y_k]+P_k(Y_1,\ldots,Y_{k-1}),
\]
where $P_1=0$ and each polynomial $P_k$ depends only on lower heights. Since $\operatorname{ad}A$ acts on $\mathfrak g_\alpha$ by the nonzero scalar $\alpha(A)$, \eqref{eq:seed-root-coordinates} is uniquely solvable one height at a time. The recursion ends at the largest positive-root height, proving existence and uniqueness of $n_A(b)$. Conjugating \eqref{eq:seed-root-coordinates} by $t\in T_\C$ gives the stated action on $b$. The root-character isomorphism \eqref{eq:root-character-isomorphism} then proves freeness and transitivity.
\end{proof}

We call $x_*(A):=x_A(1,\ldots,1)$ the \emph{reference constant reduction} for $A$, and $\widetilde f_*(w):=e^{Aw}\cdot x_*(A)$ the \emph{reference developing map} for $A$.

\begin{proof}[Proof of \cref{thm:radial-parameterization}]
\leavevmode\par
\Cref{prop:constant-reductions} shows $u^{A_\delta,x}$ is radial. Conversely, suppose that $u$ is radial. By \cref{prop:cone-adapted-logarithm}, choose an associated nondegenerate unitary integral curve $f=(\widetilde f,\rho_f)$ with monodromy $M_\delta$. By \cref{lem:radial-rotation-lift}, there is $h\in Z_K(M_\delta)$ such that
\[
    h^{-1}\cdot\widetilde f(w+\ii\theta)
      =e^{\ii\theta A_\delta}h^{-1}\cdot\widetilde f(w).
\]
Let $q(w):=e^{-A_\delta w}h^{-1}\cdot\widetilde f(w)$. This rotation equivariance gives
\[
    q(w+\ii\theta)
      =e^{-A_\delta(w+\ii\theta)}h^{-1}\cdot\widetilde f(w+\ii\theta)
      =q(w),
\]
for all $w\in\C$ and $\theta\in[0,2\pi]$. Since $q$ is holomorphic, this means $q$ is constant on $\C$, say $q\equiv x$. Since $h\in K$, the curve $f_{A_\delta,x}=h^{-1}\cdot f$ is a nondegenerate unitary integral curve. So $x$ is a constant reduction for $A_\delta$, and $u=u^{A_\delta,x}$.

Set $x_*:=x_*(A_\delta)$ and define the map in \eqref{eq:radial-parameterization} by
\[
    [t]\longmapsto u^{A_\delta,t\cdot x_*}.
\]
By \cref{prop:constant-reduction-construction,prop:constant-reductions}, this map takes values in $\mathscr R_\delta$. Since $T\subset K$ and $T$ commutes with $A_\delta$, it is well defined on $T\backslash T_\C$. The characterization proved above and the transitivity in \cref{prop:constant-reduction-construction} give surjectivity.

If $u^{A_\delta,x}=u^{A_\delta,x'}$, the Toda--flag correspondence gives $k\in K$ with $e^{A_\delta w}\cdot x'=ke^{A_\delta w}\cdot x$. Comparing \eqref{eq:constant-reduction-rotation} for the two curves and applying \cref{lem:pointwise-stabilizer} shows that $k$ commutes with $e^{\ii\theta A_\delta}$ for every $\theta$. Thus $k\in Z_K(A_\delta)=T$ by regularity of $A_\delta$, and $x'=k\cdot x$ at $w=0$. Conversely, constant $T$-translations preserve the solution. The freeness of the $T_\C$-action gives injectivity of \eqref{eq:radial-parameterization}.
\end{proof}
\section{Classification of finite-energy solutions on the punctured plane}
\label{sec:global-classification}

To prove \cref{thm:intro-classification}, we show that every solution in $\mathscr S_\delta$ is a centralizer deformation of a fixed radial solution. This deformation preserves cone exponents and total energies. The main step uses normalized algebraic lifts and the Drinfeld--Sokolov normal form.

\begin{definition}\label{def:centralizer-deformation}
Let $A\in\mathcal C_+$, and set $M_A:=e^{2\pi\ii A}$. Let $x$ be a constant reduction for $A$. For $c\in Z_G(M_A)$, let $f_{A,x;c}$ be the unitary flag curve with developing map
\begin{equation}\label{eq:centralizer-developing-curve}
    \widetilde f_{A,x;c}(w):=c e^{Aw}\cdot x,
    \qquad w\in\C.
\end{equation}
We call the associated Toda solution $u^{A,x;c}$ a \emph{centralizer deformation} of the \emph{radial seed} $u^{A,x}$.
\end{definition}

The solution $u^{A,x;c}$ is well defined. Indeed, $cM_A=M_Ac$ shows that the developing map \eqref{eq:centralizer-developing-curve} has unitary monodromy $M_A$. Constant $G$-actions preserve integrality and nondegeneracy, so the Toda--flag correspondence gives a smooth solution on $\C^*$.

\begin{proposition}\label{prop:centralizer-cone-exponents}
Let $u^{A,x;c}$ be a centralizer deformation. Then $\delta_{0,i}(u^{A,x;c})=\alpha_i(A)$ for $1\leq i\leq n$. Moreover, each of its components has finite total energy.
\end{proposition}

\begin{proof}
Put $A_i=\dd\varrho_i(A)$ and choose $v\ne0$ representing $\pi_i(x)$. On a domain with a fixed branch of $\log z$, the $i$th fundamental projective component of $f_{A,x}$ has a holomorphic lift
\[
    \vec f_i=z^{A_i}v=\sum_{j=0}^{s}z^{\lambda_j}v_j,
    \qquad \lambda_0<\cdots<\lambda_s,
\]
where $v_j\ne0$ is the projection of $v$ onto the $\lambda_j$-eigenspace of $A_i$. Since $A_i$ is Hermitian, the $\lambda_j$ are real and the $v_j$ are pairwise orthogonal. Nondegeneracy gives $s\geq1$.

For any invertible constant linear map $H$, \eqref{eq:fs-wronskian-formula} gives the local area density of $[H\vec f_i]$ as
\[
    q_H=\partial_z\partial_{\bar z}\log\|H\vec f_i\|^2
       =\frac{\|(\wedge^2H)(\vec f_i\wedge \vec f_i')\|^2}{\|H\vec f_i\|^4}.
\]
The invertibility of $H$ and $\wedge^2H$ preserves the leading orders of $\vec f_i$ and $\vec f_i\wedge \vec f_i'$, so $q_H$ and $q_I$ have the same cone exponent at $0$ as in the proof of \cref{thm:finite-area-conical}. Taking $H$ to be a linear lift of $c$ on $V_i$ and using \cref{prop:constant-reductions} proves the cone exponent equalities.

For this choice of $H$, invertibility gives constants $C_1,C_2>0$ such that $C_1q_I\leq q_H\leq C_2q_I$ on $\C^*$. Since $u^{A,x}$ has finite total energy, so does $u^{A,x;c}$.
\end{proof}

We now prove that centralizer deformations exhaust $\mathscr S_\delta$.

We begin with the construction of normalized algebraic lifts of the associated flag curves. Let $\mathcal R:=\C[z,z^{-1}]$ be the ring of Laurent polynomials in $z$, and set $D:=z\partial_z$. Under the covering map $z=e^w$, $D$ lifts to $\dd/\dd w$, denoted by the same symbol. Put $\ee:=\sum_i\ee_i$.

\begin{lemma}\label{lem:simple-root-derivative-normalization}
Let $f=(\widetilde f,\rho_f)$ be a nondegenerate unitary integral curve on $\C^*$ whose fundamental projective components all have finite total area, and whose monodromy is $M_A=e^{2\pi\ii A}$, where $A\in\ii\mathfrak t$. Then there is a holomorphic map $\widehat f:\C\to G$ satisfying
\[
    \widehat f(w)\cdot B_0=\widetilde f(w),\qquad
    \widehat f(w+2\pi\ii)=M_A\widehat f(w).
\]
The map $w\mapsto e^{-Aw}\widehat f(w)$ descends to an algebraic map on $\C^*$, and
\begin{equation}\label{eq:simple-root-normalized-derivative}
    \widehat f^{-1}D\widehat f=\ee+\beta,
    \qquad \beta\in\mathfrak b_0\otimes_\C \mathcal R.
\end{equation}
We call such a lift $\widehat f$ a \emph{normalized algebraic lift} of $f$.

Any two normalized algebraic lifts $\widehat f_1, \widehat f_2$ of $f$ are related by $\widehat f_2(w)=\widehat f_1(w)b(e^w)$ for a unique algebraic map $b:\C^*\to N_-$. Conversely, right multiplication by any such map also gives a normalized algebraic lift.
\end{lemma}

\begin{proof}
As in \cref{thm:global-flag-reduction}, the flag single-valued reduction $g$ of $f$ relative to $w\mapsto e^{Aw}$ extends holomorphically to $\Pro^1$. The extension is algebraic by Chow's theorem \cite[Appendix~C, Theorem~C.1]{EisenbudHarris3264}. Since $B_0$ is split solvable and the Picard group $\operatorname{Pic}(\C^*)=0$, the $B_0$-torsor obtained by pulling back $G\to F$ along $g$ is trivial by \cite[Proposition~16.55(a)]{MilneAlgebraicGroups}. Choosing an algebraic section of this trivial pullback torsor gives an algebraic map $\widehat g_0:\C^*\to G$ satisfying $\widehat g_0(z)\cdot B_0=g(z)$. We set $\widehat f_0(w):=e^{Aw}\widehat g_0(e^w)$.

The product rule gives
\[
    J_0:=\widehat f_0^{-1}D\widehat f_0
    =\operatorname{Ad}_{\widehat g_0^{-1}}A
      +\widehat g_0^{-1}D\widehat g_0,
\]
so $J_0$ has Laurent coefficients by algebraicity of $\widehat g_0$. Integrality and nondegeneracy give
\[
    J_0=\sum_{i=1}^n a_i\ee_i+\beta_0,
    \qquad a_i\in\mathcal R,\quad
    \beta_0\in\mathfrak b_0\otimes_\C\mathcal R,
\]
where each $a_i$ is nowhere zero on $\C^*$, hence $a_i=c_i z^{m_i}$ with $c_i\in\C^*$ and $m_i\in\mathbb Z$. By \eqref{eq:root-character-isomorphism}, there is a unique algebraic map $t:\C^*\to T_\C$ satisfying $\chi_i(t)=a_i$ for every $i$.

Set $\widehat f(w):=\widehat f_0(w)t(e^w)$. It projects to $\widetilde f$ and has monodromy $M_A$. The map $w\mapsto e^{-Aw}\widehat f(w)$ descends to the algebraic map $z\mapsto\widehat g_0(z)t(z)$ on $\C^*$. The product rule gives
\[
    \widehat f^{-1}D\widehat f
       =\operatorname{Ad}_{t^{-1}}J_0+t^{-1}Dt
       =\ee+\beta,
    \qquad \beta\in\mathfrak b_0\otimes_\C\mathcal R.
\]

For the remaining freedom, let $\widehat f_1,\widehat f_2$ be two normalized algebraic lifts of $f$. Write $\widehat f_j(w)=e^{Aw}\widehat g_j(e^w)$, where $\widehat g_j:\C^*\to G$ are algebraic. Since both $\widehat g_j$ lift $g$, we have $\widehat g_1(z)B_0=\widehat g_2(z)B_0$. Hence $b:=\widehat g_1^{-1}\widehat g_2$ is an algebraic $B_0$-valued map and satisfies $\widehat f_2(w)=\widehat f_1(w)b(e^w)$. This identity determines $b$ uniquely.  Let $t_b$ be its image in $T_\C$. Since $N_-$ acts trivially on $(\mathfrak b_0+\sum_i\C\ee_i)/\mathfrak b_0$, the product rule gives
\[
    \widehat f_2^{-1}D\widehat f_2
    =\operatorname{Ad}_{b^{-1}}
       (\widehat f_1^{-1}D\widehat f_1)+b^{-1}Db
    \equiv\sum_i\chi_i(t_b)^{-1}\ee_i
       \pmod{\mathfrak b_0}.
\]
Normalization and \eqref{eq:root-character-isomorphism} therefore force $t_b=1$, so $b$ takes values in $N_-$. Conversely, an algebraic $N_-$-valued right factor preserves the projection, monodromy and algebraicity after removal of $e^{Aw}$, while the same calculation preserves normalization.
\end{proof}

By \cref{lem:simple-root-derivative-normalization}, normalized algebraic lifts differ by right multiplication by $N_-$-valued algebraic maps. We use this freedom to put the logarithmic derivative in Drinfeld--Sokolov normal form, and then prove that its coefficients are constant. This yields an explicit expression for the developing map in terms of constant data.

Write $2\rho^\vee=\sum_i r_i\hh_i$ and set $\mathsf p:=\sum_i r_i\ff_i$, so that $(\ee,2\rho^\vee,\mathsf p)$ is the principal $\mathfrak{sl}_2$-triple. Set $\mathfrak s:=\ker(\operatorname{ad}\mathsf p)$ and fix a basis $s_1,\ldots,s_n$ of $\mathfrak s$; see \cite[Section~4.4, Theorem~5 and Lemma~12]{Kostant1963}. The following lemma records the holomorphic form of the Drinfeld--Sokolov normal-form theorem \cite[Proposition~6.1]{DrinfeldSokolov} in these conventions.

\begin{lemma}[Drinfeld--Sokolov normal form]
\label{lem:DS-normal-form}
Let $\beta$ be a holomorphic $\mathfrak b_0$-valued function on a domain $\Omega$ in the $w$-plane. There is a unique holomorphic map $b_{\mathrm{DS}}:\Omega\to N_-$ such that
\begin{equation}\label{eq:DS-normal-form}
    b_{\mathrm{DS}}^{-1}(D+\ee+\beta)b_{\mathrm{DS}}
       =D+\ee+\sum_{a=1}^n\mathcal W_as_a.
\end{equation}
The components of $\log b_{\mathrm{DS}}$ and the coefficients $\mathcal W_a$ are differential polynomials in $\beta$ determined by these $s_a$. These $\mathcal W_a$ are invariant under $N_-$ gauge transformations.

In the special case $\beta=\mu:\Omega\to\mathfrak h$, these $\mathcal W_a$ have the form
\begin{equation}\label{eq:Miura-differential-polynomials}
    \mathcal W_a=\mathcal W_a(\mu,D\mu,\ldots,D^{N_a}\mu),
\end{equation}
where $N_a$ is the highest order of derivative of $\mu$ appearing in $\mathcal W_a$.
\end{lemma}

For Cartan-valued input, this is the normal form in \cite[Lemma~4.1(i)--(ii)]{KarmakarLinNieWei}, written in our negative Borel convention. The coefficients $\mathcal W_a$ are called the \emph{$W$-invariants} there.

\begin{proposition}\label{prop:constant-flag-generator}
Under the hypotheses of \cref{lem:simple-root-derivative-normalization}, choose a normalized algebraic lift $\widehat f$. There exist an algebraic map $b_{\mathrm{DS}}:\C^*\to N_-$ and constants $c_a\in\C$ such that
\[
    \widehat f_{\mathrm{can}}(w):=\widehat f(w)b_{\mathrm{DS}}(e^w), \qquad
    \widehat f_{\mathrm{can}}^{-1}D\widehat f_{\mathrm{can}}
    =\Lambda:=\ee+\sum_{a=1}^n c_a s_a.
\]
Moreover, there exist $d\in G$ and $A'\in\mathcal C_+$ such that
\begin{equation}\label{eq:constant-flag-generator}
    \widetilde f(w)=de^{A'w}d^{-1}\cdot\widetilde f(0).
\end{equation}
\end{proposition}

\begin{proof}
Applying \cref{lem:DS-normal-form} to \eqref{eq:simple-root-normalized-derivative}, we obtain
\[
    \widehat f_{\mathrm{can}}^{-1}D\widehat f_{\mathrm{can}}
      =\ee+\sum_{a=1}^n\mathcal W_as_a.
\]
Because $\beta\in\mathfrak b_0\otimes_\C\mathcal R$ and $D(\mathcal R)\subset\mathcal R$, we have $b_{\mathrm{DS}}\in N_-(\mathcal R)$ and $\mathcal W_a\in\mathcal R$. We now prove that the $\mathcal W_a$ are constant.

Write $\widehat f(w)=e^{Aw}\widehat g(e^w)$ with $\widehat g\in G(\mathcal R)$. In a weight basis of the faithful adjoint representation, $\operatorname{Ad}_{e^{Aw}}$ is diagonal with entries $e^{\lambda(A)w}$, where $\lambda$ ranges over the weights. Since $A\in\ii\mathfrak t$, all $\lambda(A)$ are real. By algebraicity of $\widehat g$, the entries of $\operatorname{Ad}_{\widehat g(e^w)}$ are Laurent polynomials in $e^w$. Hence the entries of $\operatorname{Ad}_{\widehat f(w)}$ are finite sums of exponentials $e^{\eta w}$ with $\eta\in\R$.

Set $h(w):=\widehat f(0)^{-1}\widehat f(w)$. Since $h(0)=1$, we have the Gauss decomposition $h=n_+tn_-$ near $w=0$, with $n_\pm\in N_\pm$ and $t\in T_\C$ \cite[Theorem~21.84]{MilneAlgebraicGroups}. The factors are rational in the matrix coefficients of $\operatorname{Ad}_h$. Put $h_{\mathrm M}:=hn_-^{-1}=n_+t\in T_\C\ltimes N_+$. Right multiplication by $n_-^{-1}$ preserves normalization, so $h_{\mathrm M}^{-1}Dh_{\mathrm M}=\ee+\mu$ with $\mu\in\mathfrak h$. Projection to $(T_\C\ltimes N_+)/N_+\cong T_\C$ gives
\begin{equation}\label{eq:Gauss-Cartan-logarithmic-derivative}
    \mu=t^{-1}Dt,\qquad
    \alpha_i(\mu)=\frac{Dy_i}{y_i},
    \qquad y_i:=\chi_i(t).
\end{equation}
By rationality of the factors, each $y_i$ is a rational combination of exponentials $e^{\eta w}$ with $\eta\in\R$, and hence meromorphic in $\C$. By $N_-$ gauge invariance, these $\mathcal W_a$ satisfy \eqref{eq:Miura-differential-polynomials} with this $\mu$.

Restrict to $w=x\in\R$. Dominance of the extreme terms shows that $y_i$ has no zeros or poles for sufficiently large $|x|$. So $\alpha_i(\mu)=Dy_i/y_i$ has finite limits as $x\to\pm\infty$, and all its positive-order derivatives tend to zero. Since the simple roots form a basis of $\mathfrak h^*$, $\mu$ and the derivatives occurring in \eqref{eq:Miura-differential-polynomials} are bounded at both ends. Hence $\mathcal W_a(e^x)$ is bounded as $x\to\pm\infty$. Since $\mathcal W_a\in\mathcal R$, these bounds force $\mathcal W_a$ to be a constant $c_a$.

It remains to prove \eqref{eq:constant-flag-generator}. Since $\widehat g b_{\mathrm{DS}}\in G(\mathcal R)$, the matrix entries of $\Gamma(w):=\operatorname{Ad}_{\widehat f_{\mathrm{can}}(w)}$ are finite sums of exponentials with real exponents, by the same argument. Let $S\subset\R$ be their finite set of exponents and put $P(\xi):=\prod_{\eta\in S}(\xi-\eta)$. Integrating the logarithmic derivative gives $\widehat f_{\mathrm{can}}(w)=\widehat f_{\mathrm{can}}(0)e^{w\Lambda}$, so $D\Gamma=\Gamma\operatorname{ad}\Lambda$, and hence
\[
    0=P(D)\Gamma=\Gamma\,P(\operatorname{ad}\Lambda).
\]
Invertibility of $\Gamma$ gives $P(\operatorname{ad}\Lambda)=0$. The polynomial $P$ has distinct real roots, so $\Lambda$ is semisimple with real adjoint spectrum. By \cite[Section~4.2, Lemma~10]{Kostant1963}, $\Lambda\in\ee+\mathfrak s\subset\ee+\mathfrak b_0$ is regular. Hence we may write $\Lambda=\operatorname{Ad}_qA'$ for some $q\in G$ and $A'\in\mathcal C_+$. Projecting $\widehat f_{\mathrm{can}}$ to $F$ gives \eqref{eq:constant-flag-generator}, with $d:=\widehat f_{\mathrm{can}}(0)q$.
\end{proof}

\begin{proof}[Proof of \cref{thm:intro-classification} part~(1)]
\leavevmode\par
Set $x_*:=x_*(A_\delta)$, as defined after \cref{prop:constant-reduction-construction}. By \cref{prop:centralizer-cone-exponents}, $u^{A_\delta,x_*;c}\in\mathscr S_\delta$ for every $c\in Z_G(M_\delta)$. Since constant $Z_K(M_\delta)$-translations preserve the associated Toda solution, the following map is well defined:
\begin{equation}\label{eq:finite-energy-classification}
    \Phi_\delta:Z_K(M_\delta)\backslash Z_G(M_\delta)\longrightarrow\mathscr S_\delta,
    \qquad [c]\longmapsto u^{A_\delta,x_*;c}.
\end{equation}

To prove surjectivity, let $u\in\mathscr S_\delta$. By \cref{prop:cone-adapted-logarithm}, choose an associated nondegenerate unitary integral curve $f=(\widetilde f,\rho_f)$ with monodromy $M_\delta$. Applying \cref{prop:constant-flag-generator}, write $\widetilde f$ as in \eqref{eq:constant-flag-generator}. Since the $G$-action preserves integrality and nondegeneracy, $d^{-1}\cdot\widetilde f(0)$ is a constant reduction for $A'$. By \cref{prop:constant-reductions}, the associated Toda solution has cone exponents $\alpha_i(A')$ at $0$. The local calculation in the proof of \cref{prop:centralizer-cone-exponents}, valid for any invertible constant linear map, gives $\alpha_i(A')=\delta_{0,i}(u)=\delta_i$ for every $i$. Hence $A'=A_\delta$.

Set $\widetilde f_*(w):=e^{A_\delta w}\cdot x_*$. By \cref{prop:constant-reduction-construction}, $d^{-1}\cdot \widetilde f(0)=t\cdot x_*$ for some $t\in T_\C$. Thus $\widetilde f=c\cdot\widetilde f_*$ with $c=dt$. Equivariance implies that $(cM_\delta)^{-1}M_\delta c$ fixes $\widetilde f_*$ pointwise. By \cref{lem:pointwise-stabilizer}, it is the identity, so $c\in Z_G(M_\delta)$. Hence $u=\Phi_\delta([c])$, proving surjectivity.

Suppose $u^{A_\delta,x_*;c}=u^{A_\delta,x_*;c'}$. The Toda--flag correspondence gives $h\in K$ with $c'\cdot\widetilde f_*=hc\cdot\widetilde f_*$. Since both $c\cdot\widetilde f_*$ and $c'\cdot\widetilde f_*$ have monodromy $M_\delta$, its uniqueness gives $h\in Z_K(M_\delta)$. The element $c'^{-1}hc$ fixes $\widetilde f_*$ pointwise, so \cref{lem:pointwise-stabilizer} gives $c'=hc$, proving injectivity of \eqref{eq:finite-energy-classification}.

Every $t\in T_\C\subset Z_G(M_\delta)$ commutes with $A_\delta$, so $u^{A_\delta,x_*;t}=u^{A_\delta,t\cdot x_*}$. By \cref{thm:radial-parameterization}, $\Phi_\delta$ therefore restricts to a bijection from $T\backslash T_\C$ onto $\mathscr R_\delta$.
\end{proof}

\begin{proof}[Proof of \cref{thm:intro-classification} part~(2)]
\leavevmode\par
For $u\in\mathscr S_\delta$, by part~(1), choose $c\in Z_G(M_\delta)$ such that an associated developing map is
\[
    \widetilde f(w)=ce^{A_\delta w}\cdot x_*.
\]
In the coordinate $\zeta=1/z=e^v$ at infinity, the developing map is $ce^{-A_\delta v}\cdot x_*$. The same argument as in the proof of \cref{prop:centralizer-cone-exponents} gives $\delta_{\infty,i}(u)=\alpha_i(-\kappa A_\delta)$, since $-\kappa A_\delta$ is the unique element of $\mathcal C_+$ that is $K$-conjugate to $-A_\delta$.

The estimate \eqref{eq:cone-radial-derivative}, applied also to $u_{i,\infty}(\zeta)=u_i(1/\zeta)-4\log|\zeta|$, gives the following limits uniformly in $\arg z$:
\begin{align*}
    r\partial_r u_i&\longrightarrow2\bigl(\delta_{0,i}(u)-1\bigr) &&(r\to0),\\
    r\partial_r u_i&\longrightarrow-2\bigl(\delta_{\infty,i}(u)+1\bigr) &&(r\to\infty).
\end{align*}
Each component has finite total energy because it has conical singularities at both ends. Integrating \eqref{eq:toda-local} on annuli and applying Green's formula yields the identity
\begin{equation}\label{eq:total-energy}
    \sum_j c_{ij}\mathrm{En}_j(u)
      =\pi\bigl(\delta_{0,i}(u)+\delta_{\infty,i}(u)\bigr).
\end{equation}
Using $\alpha_i=\sum_jc_{ij}\omega_j$ and invertibility of the Cartan matrix $\mathsf C$, we obtain the asserted total energy formula.
\end{proof}

\begin{remark}\label{rem:comparison-KLNW}
Our approach is inspired by \cite{KarmakarLinNieWei}, especially its use of Drinfeld--Sokolov normal form. The common strategy is to use this normal form to match a general solution with a model determined by the singularities.

As in \cite[(4.8)--(4.11)]{KarmakarLinNieWei}, we use a torus factor to set the simple-root coefficients of a logarithmic derivative equal to one, working here with an algebraic lift of the developing flag curve.

Both proofs obtain constant Drinfeld--Sokolov data from bounds at the two ends. In \cite[Proposition~5.1]{KarmakarLinNieWei}, finite-energy estimates for derivatives of $U_i$ bound $z^{d_j}W_j$. We bound the Laurent polynomials $\mathcal W_a$ using the logarithmic derivatives of the torus factor in the Gauss decomposition.

Karmakar--Lin--Nie--Wei compute the Drinfeld--Sokolov coefficients from the singularity data at the origin and match them with those of the prescribed model \cite[Theorems~5.3 and~6.1]{KarmakarLinNieWei}. We identify the model by comparing its cone exponents with those of the original solution.

The formulas for total energies and cone exponents at infinity agree with \cite[Theorems~1.1 and~1.2]{KarmakarLinNieWei} under the following notation correspondence:
\[
\begin{gathered}
a_{ij}^{\mathrm{KLNW}}=c_{ij},\qquad
\gamma_i^{\mathrm{KLNW}}=\delta_i-1,\qquad
\mu_i^{\mathrm{KLNW}}=\delta_i, \\
w_0^{\mathrm{KLNW}}=A_\delta,\qquad
E_0^{\mathrm{KLNW}}=\rho^\vee,\qquad
\sigma_i^{\mathrm{KLNW}}=\frac{2}{\pi}\mathrm{En}_i.
\end{gathered}
\]
Our Borel convention is opposite to theirs. Their fundamental representations have highest weights $\omega_i$, whereas ours have lowest weight $-\omega_i$ and hence highest weight $-\kappa\omega_i$.
\end{remark}

As a final consequence of \cref{thm:intro-classification}, the dimension of the parameter space can be read off from the cone exponents $\delta_i$ at the origin by counting the positive roots $\alpha=\sum_i m_i\alpha_i$ for which $\sum_i m_i\delta_i$ is an integer.

\begin{corollary}\label{cor:radiality-from-cone-exponents}
Fix $\delta\in(\R_{>0})^n$ and set $\Delta^+_\delta:=\{\alpha\in\Delta^+:\alpha(A_\delta)\in\mathbb Z\}$ and $N_\delta:=|\Delta^+_\delta|$. The parameter spaces for $\mathscr S_\delta$ and $\mathscr R_\delta$ have real dimensions $n+2N_\delta$ and $n$, respectively. In particular, $\mathscr S_\delta=\mathscr R_\delta$ if and only if $N_\delta=0$.
\end{corollary}

\begin{proof}
Let $\mathfrak p_\delta$ be the subspace of $\ii\mathfrak k$ fixed by $\operatorname{Ad}_{M_\delta}$. Since $M_\delta\in K$, the Cartan decomposition $K\times\ii\mathfrak k\xrightarrow{\sim}G,\ (k,X)\mapsto k\exp X$ is equivariant under conjugation by $M_\delta$. It therefore restricts to $Z_K(M_\delta)\times\mathfrak p_\delta\xrightarrow{\sim}Z_G(M_\delta)$. Passing to the left quotient gives $\mathfrak p_\delta\xrightarrow{\sim}Z_K(M_\delta)\backslash Z_G(M_\delta)$ by $X\mapsto[\exp X]$. This bijection sends $\ii\mathfrak t$, which is $n$-dimensional, onto $T\backslash T_\C$.

Since $\mathfrak g=\mathfrak k\oplus\ii\mathfrak k$ as real vector spaces, the complexification of $\mathfrak p_\delta$ is the subalgebra of $\mathfrak g$ fixed by $\operatorname{Ad}_{M_\delta}$. Since $\operatorname{Ad}_{M_\delta}$ acts on $\mathfrak g_\alpha$ by multiplication by $e^{2\pi\ii\alpha(A_\delta)}$, we have
\[
    \mathfrak p_\delta\otimes_\R\C
    =\mathfrak h\oplus
      \bigoplus_{\alpha\in\Delta^+_\delta}
      \bigl(\mathfrak g_\alpha\oplus\mathfrak g_{-\alpha}\bigr).
\]
Each root space has complex dimension one, so $\dim_\R\mathfrak p_\delta=n+2N_\delta$.
\end{proof}

\begin{remark}
In type $A_n$, the positive roots are $\alpha=\alpha_i+\cdots+\alpha_j$ for some $1\leq i\leq j\leq n$, so $\alpha(A_\delta)=\delta_i+\cdots+\delta_j$. Thus $N_\delta=0$ exactly when none of the consecutive sums $\delta_i+\cdots+\delta_j$ is an integer. With $\gamma_k=\delta_k-1$, the implication from $N_\delta=0$ to radiality of every solution recovers the criterion of \cite[Corollary~1.2]{LinWeiYe}.
\end{remark}

\paragraph{Acknowledgements.}
GPT-6 Astra supplied the core argument for the finite-energy propagation proof in \cref{sec:energy-propagation}. GPT-6.1 Sol helped identify the Drinfeld--Sokolov normal forms used in \cite{KarmakarLinNieWei} as a tool for the global classification in \cref{sec:global-classification} and suggested ways to apply them. GPT-5.6 Sol and GPT-6.1 Sol assisted with drafting and revising the manuscript. The authors take full responsibility for the mathematical content and the final text.

Yiqian Shi is supported in part by the National Natural Science Foundation of China (NSFC, Grant No.~11931009). Bin Xu is supported in part by the Project of Stable Support for Youth Team in Basic Research Field of the Chinese Academy of Sciences (Grant No.~YSBR-001) and by NSFC (Grant No.~12271495).

\end{document}